\documentclass{article}
\usepackage{graphicx} 
\usepackage[dvipsnames]{xcolor}
\usepackage[utf8]{inputenc}

\usepackage{todonotes}

\usepackage{amsmath,amsfonts,amssymb,amsthm}
\usepackage{url}
\usepackage{tikz}
\usetikzlibrary{calc}

\newcommand{\tikzmark}[1]{\tikz[overlay,remember picture] \node (#1) {};}

\title{Perspectives on locally weighted ensemble Kalman methods}
\author{Philipp Wacker\thanks{philipp.wacker@canterbury.ac.nz}\\School of Mathematics and Statistics\\University of Canterbury\\
\url{https://orcid.org/0000-0001-8718-4313} }
\date{\today}

\DeclareMathOperator{\argmin}{argmin}
\newcommand{\R}{\mathbb R}
\newcommand{\N}{\mathbb N}
\newcommand{\X}{\mathcal X}
\newcommand{\Y}{\mathcal Y}
\newcommand{\eps}{\varepsilon}
\renewcommand{\d}{\mathrm d}

\newcommand{\new}[1]{{#1}}
\newcommand{\neww}[1]{{#1}}

\newtheorem{lem}{Lemma}
\newtheorem{definition}{Definition}
\newtheorem{rem}{Remark}
\usepackage[
backend=biber,
style=authoryear 
]{biblatex}
\begin{document}

\maketitle

\begin{abstract}
    This manuscript derives locally weighted ensemble Kalman methods from the point of view of ensemble-based function approximation. This is done by using pointwise evaluations to build up a local linear or quadratic approximation of a function, tapering off the effect of distant particles via local weighting. This introduces a candidate method (the locally weighted Ensemble Kalman method for inversion) with the motivation of combining some of the strengths of the particle filter (ability to cope with nonlinear maps and non-Gaussian distributions) and the Ensemble Kalman filter (no filter degeneracy). \neww{We provide some numerical evidence for the accuracy of locally weighted ensemble methods, both in terms of approximation and inversion.}
\end{abstract}

\textit{Keywords: Ensemble Kalman, Filtering, Inversion, Optimisation, Sampling, Bayesian Inference}

\textit{MSC class: 62F15, 65N75, 37N40, 90C56}

\tableofcontents

\section{Introduction}
The following problem motivates the study of Ensemble Kalman methods: Let $\X$ be a Hilbert space, the parameter space and a possibly nonlinear measurement operator $A: \X\to \Y$ into a feature (or observation) space $\Y$. Given an observation $y\in \Y$ from a noisy measurement
\begin{equation}
    y = A(u) + \eps,
\end{equation}
where $\eps\in \Y$ is a realisation of measurement noise (often assumed to be Gaussian), we want to infer the original parameter $u\in \X$. Often, a misfit functional $\Phi: \X \to \R$ is constructed which associates a \textit{cost} to a given parameter value. For example, when $\eps\sim \mathcal N(0,I)$ is a Gaussian random variable (setting covariance to the identity purely out of convenience),
\begin{equation}
    \Phi(u) = \frac{1}{2}\|y - A(u)\|^2
\end{equation}

This setting is quite general and covers many interesting optimisation, filtering, and sampling problems. For concreteness and in order to exclude a few pathologies of infinite-dimensional spaces (although the Ensemble Kalman method is quite robust in this regard), we will set $\X = \R^{p}$ and $\Y = \R^d$.

In this manuscript we will consider in a very general sense the \textit{family of Ensemble Kalman methods}, which have been developed in the context of filtering first (\cite{evensen2003ensemble}), based on the Kalman filter (EnKF, \cite{kalman}), and then applied to inversion/optimisation (EKI, \cite{iglesias2013ensemble}), and sampling (EKS, \cite{garbuno2020interacting} and ALDI, \cite{garbuno2020affine}), with applications in weather forecasting, geophysics, control, robotics, machine learning, any many others. We refer the reader to \cite{calvello2022ensemble} for a recent overview and the references therein. We will not attempt to enumerate all relevant or interesting variants which have been established, like the original EnKF, the Ensemble Square Root filter (EnSRF, \cite{tippett2003ensemble}), the EKS, the EKI, ALDI, Tikhonov-regularised Ensemble Kalman inversion (TEKI, \cite{chada2020tikhonov}), and countless others, with both discrete-time and continuous time variants, in the deterministic and stochastic or optimal transport setting. Instead, we focus our attention on the Ensemble Kalman method for inversion, the EKI, as a test bed. 

In this manuscript, when we talk about ``the (linear) EKI'', we will mean\footnote{accepting that this falls short of accommodating all variants of the Ensemble Kalman method family} the following coupled system of $J\in \N$ ordinary differential equations (ODEs). 
\begin{equation}\label{eq:EK_evolution}
    \frac{\d}{\d t} u^{(i)} = -C_{u,A} \cdot \left(A(u^{(i)}) - y\right).
\end{equation}
This system of equations tracks the evolution of an \textit{ensemble} of particles $\mathcal E = \{u^{(i)}(t)\}_{i=1}^J$, coupled via the empirical parameter-feature covariance matrix
\begin{align*}
    C_{u,A} &= \frac{1}{J-1}\sum_{j=1}^J \left(u^{(j)} - \mu\right) \cdot  \left(A(u^{(j)})  - \mu_A\right)^{\neww\top}
\end{align*}
and we further define
\begin{alignat*}{2}
    &\text{empirical mean } &\quad\mu &= \frac{1}{J}\sum_{j=1}^J u^{(j)},\\
    &\text{empirical feature mean} &\quad \mu_A &= \frac{1}{J}\sum_{j=1}^J A(u^{(j)}), \text{ and }\\
    &\text{empirical covariance}&C_{u,u} &= \frac{1}{J-1}\sum_{j=1}^J \left(u^{(j)} - \mu\right) \cdot  \left(u^{(j)}  - \mu\right)^\top.
\end{alignat*}
Note that, of course, most quantities here depend on time via their dependency on the time-varying ensemble $\{u^{(i)}(t)\}_{i=1}^J$, although we suppress this dependence notationally.

This manuscript considers a variation of the EKI in \eqref{eq:EK_evolution} in the following form:
\begin{equation}\label{eq:EK_evolution_loc}
    \boxed{\frac{\d}{\d t} u^{(i)} = -C_{u,A}^\kappa{\color{red}(u^{(i)})} \cdot \left(A(u^{(i)}) - y\right)}
\end{equation}

Here, 
\begin{align}
    C_{u,A}^\kappa{\color{red}(x)} &:= \sum_{j=1}^J \kappa^{(j)}{\color{red}(x)}\left(u^{(j)} - \mu^\kappa{\color{red}(x)}\right) \cdot  \left(A(u^{(j)})  - \mu_A^\kappa{\color{red}(x)}\right)^\top
\end{align}
is a  $\kappa$-weighted empirical ensemble covariance, with weights anchored to the reference point ${\color{red}x}$$\kappa^{(j)}{\color{red}(x)}$, depending on the distance between ensemble members and $x$ (where the dependence on the reference point is highlighted in {\color{red}red}). Alternatively, weights can be defined in some other way without using a distance. $\kappa$ can be interpreted as a $J$-simplex probability distribution field (i.e. $\X\mapsto \mathcal P(\{1,\ldots, J\})$). Similarly, $\mu^\kappa(x)$ is an empirical mean of the ensemble, locally weighted by distance to the reference point. By taking the large ensemble limit (or by replacing the empirical measure given by the ensemble with a more general probability measure), we can also directly introduce mean \textit{fields} and covariance \textit{tensor fields}. This allows us to lift the method to the level of a mean field method, in the spirit of \cite{calvello2022ensemble}. A detailed explanation of all terms involved follow below, and we start by motivating why such a local weighting procedure is worthy of investigation.

The contributions of this manuscript are as follows:
\begin{itemize}
    \item Section \ref{sec:EK} gives some reasons for why the EKI is so successful in practice, motivating the sections to follow. We continue to point out some interesting connections to finite frame theory and Riemannian geometry.
    \item The main section \ref{sec:linearapprox} defines and analyses the object 
    \[ D_\kappa  A(x)= C_{A,u}^\kappa(x)\circ C_{u,u}^\kappa(x)^{-1},\]
    which is an ensemble-based gradient approximation controlled by the weight function $\kappa$. \new{By choosing a suitable weight function $\kappa$, we can perform gradient approximation from pointwise evaluation beyond the linear setting.} In particular, a uniform weight function recovers the ``statistical linearisation'' $C_{A,u}C_{u,u}^{-1}$, while the case of vanishing kernel bandwidth recovers the true gradient in the mean-field limit. We also show how these ideas can be used to construct higher-order gradient approximations, explicitly shown for the second derivative of $A$, \new{hinting at possible generalisations to higher-order optimisation and sampling methods. We perform some numerical experiments to test gradient approximation quality.}
    \item Section \ref{sec:lwEK} uses this idea to define a locally weighted form of the EKI and tests its performance on some model problems in inversion and filtering. \new{Building up on the analysis of $D_\kappa A$ in section \ref{sec:linearapprox}, we will see in lemma \ref{lem:precond_gradflow_locweight} that \textbf{locally weighted EKI performs a preconditioned gradient flow similar to ``vanilla'' EKI}.  This means the locally weighted EKI indeed generalises the EKI in a very natural way and retains its most important properties.}
\end{itemize}

\section{Ensemble Kalman-type methods}\label{sec:EK}
\subsection{Gradient descent via empirical covariance}
The power and popularity of the EKI originates in the amazing fact that \textit{in the linear setting}, where the mapping $A$ is linear (and thus can be represented by a matrix $A\in \R^{d\times p}$), \eqref{eq:EK_evolution} can be seen (see \cite{schillings2017analysis}) to be equal to the following representation:
\begin{equation}\label{eq:EK_evolution_gradient}
    \frac{\d}{\d t} u^{(i)} = -C_{u,u} \cdot \nabla \Phi(u^{(i)}) = -C_{u,u}A^\top(Au^{(i)} - y).
\end{equation}
This means the EKI performs particle-wise gradient descent on the misfit functional $\Phi$, preconditioned by the empirical covariance matrix $C_{u,u}$ of the ensemble in parameter space. The empirical parameter-feature covariance $C_{u,A}$ mediates the linear relationship between parameters and features so  $C_{u,A}\cdot \left(A(u^{(i)}) - y\right)$ indeed characterises a (preconditioned) gradient of the misfit functional.

This is surprising at first glance because \eqref{eq:EK_evolution} only requires pointwise evaluation of $A(u^{(i)}(t))$ and no explicit gradients of $\Phi$. A little reflection of course shows this is due to the fact that the gradient of $\Phi$ has a very simple form: $\nabla \Phi(u) = A^\top(A(u)-y)$.

\subsection{Ensemble coherence versus nonlinearity}
It seems natural to ask whether this can be leveraged for nonlinear forward maps $A$ in the same way. Unfortunately, this is not generally the case: Let $A:\R\to \R$ with $A(u) = u^2$, and consider an ensemble $\{u^{(i)}\}_{i=1}^J$ distributed (sufficiently) symmetrically around $0$. Then $C_{u,A}$ will be (close to) $0$, since the best-fit linear regression mapping $u\mapsto mu+b$ through the dataset $\{(u^{(i)}, A(u^{(i)}))\}_{i=1}^J$ has slope $m=0$. This means at this point (for this configuration of ensemble), the EKI terminates regardless of the value of $y$ and does not, in fact, perform a gradient descent. This is the issue with combining (strongly) nonlinear forward mappings with the EKI. See \cite{ernst2015analysis} for a more in-depth exploration of this issue.

On the other hand, by Taylor's theorem, every sufficiently smooth nonlinear function can be approximated \textit{locally} by a linear mapping. This means if the ensemble is sufficiently concentrated in a region where linearity holds approximately, the EKI behaves in a similar way as with for linear forward maps. This is not too surprising, since the EKI only ``sees'' the mapping $A$ through evaluations  $A(u^{(i)})$, so if $A$ is only evaluated in a region where $A$ is locally ``linear enough'', \eqref{eq:EK_evolution_gradient} still holds approximately.

This trade-off can be summarised as follows:
\begin{itemize}
    \item A widely distributed ensemble is able to explore the state space better, leads to quicker convergence of the EKI (due to the preconditioning via the ensemble covariance), and is sometimes unavoidable if the initial ensemble is drawn from a wide prior measure.
    \item A strongly contracted ensemble allows to approximate a given nonlinear forward mapping $A$ much better by its local linear approximation as given by its affine-linear regression mapping through $\{(u^{(i)}, A(u^{(i)}))\}_{i=1}^J$ and thus does not (completely) break the equivalence between \eqref{eq:EK_evolution} and \eqref{eq:EK_evolution_gradient}.
\end{itemize}

Another aspect of the EKI's strength is its \textit{cohesion}: The particle filter (see \cite{van2019particle}) also handles an ensemble of particles, but they are treated independently of each other, weighting them with a likelihood term. This can easily lead to filter degeneracy because of a majority of particles getting negligible weights. While there are ways of mitigating this effect, like resampling and rejuvenation, filter degeneracy is widely considered still a serious problem. The EKI on the other hand does not suffer from this problem at all: Particles are coupled together via their covariance, and the ensemble tends to contract during its time evolution (see \cite{schillings2017analysis,blomker2019well}, which pushes particles towards regions of higher probability (or lower cost function value, respectively). On the other hand, the particle filter can handle nonlinear and non-Gaussian measures, while the EKI struggles with these settings, exactly because its ensemble coherence implicitly enforces a globally linear and Gaussian structure on the problem. In other words, the EKI's cohesive property is both its strength and its weakness, just like the particle filter's lack of cohesion is both its strength and its weakness. 

The idea at the core of this manuscript (which was already introduced in \cite{reich2021fokker} in the context of sampling, and in \cite{wagner2022ensemble} for the Ensemble Kalman method in the context of rare event estimation) is to combine the strengths of wide and contracted ensembles by allowing for wide ensembles, but tapering down the effect of particles on each other via a (possibly distance-based) kernel function in \eqref{eq:EK_evolution}. In a way, this combines the strengths of the EKI and the particle filter, but in a way which is very different to \cite{frei2013bridging} and \cite{stordal2011bridging}.

We define a class of suitable kernel functions next.
\begin{definition}\label{def:kernel}
    Let $k:\X \times \X \to \R$ be a kernel function with radial structure $k(u,v) = K(\|u-v\|)$ for a symmetrical function $K:\R\to [0,\infty)$. This creates a system of locally weighting particles from a reference position $u$:
\begin{equation}\label{eq:weights}
    \kappa^{(i)}(x) = \frac{k(x,u^{(i)})}{\sum_{j=1}^J k(x,u^{(j)})}
\end{equation}
Note that $\sum_i \kappa^{(i)}(x) = 1$ and we recover uniform weights $1/J$ from the flat kernel $k(u,v) = 1$.

Sometimes we will consider a family of kernels parameterised by a bandwidth $r$. In this case we assume
\begin{enumerate}
    \item $k_r$ is a kernel function in the sense above for each $r> 0$.
    \item For every $r>0$, we have $\int k_r(x,z) dz= C$ for some constant $C$ independent of $r$ and $x$.
    \item $\lim_{r\searrow 0}  \int \|x-z\|^2 k_r(x,z) \d \nu(z) = 0$ for any fixed absolutely continuous probability measure $\nu$.
\end{enumerate}
\end{definition}

This kernelisation idea can be leveraged to write down a new coupled system of ODEs
\begin{equation}\label{eq:lwEK}
    \frac{\d}{\d t} u^{(i)} = -C_{u,A}^\kappa(u^{(i)})\cdot \left(A(u^{(i)}) - y\right),
\end{equation}
where the following quantities are locally weighted versions of the usual empirical mean and covariance terms:
\begin{equation}\label{eq:locally_weighted_moments}
\begin{split}
    C_{u,A}^\kappa(x) &= \sum_{j=1}^J \kappa^{(j)}(x) (u^{(j)} - \mu^\kappa(x)) \cdot (A(u^{(j)}) - \mu_A^\kappa(x))^\top\\
    \mu^\kappa(x) &= \sum_{j=1}^J \kappa^{(j)}(x) u^{(j)}\\
    \mu_A^\kappa(x) &= \sum_{j=1}^J \kappa^{(j)}(x) A(u^{(j)})
\end{split}
\end{equation}
It can be helpful to understand these as empirical properties of the ensemble, as seen from a reference point $x$ through a thick fog. For example, $\mu^\kappa(x)$ is the center of mass of the ensemble weighted by the particles' distance to the reference point $x$ (discounting the effect of particles far away from $x$ because they are not ``visible through the fog as seen from the point of view of $x$''). Essentially, $\mu^\kappa, \mu_A^\kappa, C_{u,u}^\kappa$ and $C_{u,A}^\kappa$ are now scalar fields and tensor fields (instead of a global mean and covariance tensor). For example, $C_{u,u}^\kappa: \X \to \mathcal L(\X\times \X, \R)$ is a $(2,0)$ tensor field (in physics notation). 

In many of the following proofs we will apply the following equality:
\begin{equation}
    \label{eq:fundamental}
    \sum_{j}\kappa^{(j)}(x)(u^{(j)} - \mu^\kappa(x)) = 0,
\end{equation}
which follows from the definition of $\mu^\kappa$ and the fact that the coefficients $\kappa^{(j)}$ sum to $1$.

\subsection{Finite frames given by an ensemble}
Before we start analysing locally weighted EK, and the quantities in \eqref{eq:locally_weighted_moments}, we recall some notions from finite frame theory (\cite{casazza2013introduction}) which will provide some arguments later in the manuscript. We will generally assume that $J > \dim(\X)$ so the ensemble is able to affinely span $\X$. It is possible to restrict $\X$ to this affine span, anyway, since the EKI does not leave this affine subspace, see the discussion in \cite{iglesias2013ensemble}. Also, we exclude the (probability $0$ for any randomly sampled initial ensemble) possibility that the affine span does not equal full space $\X$). Regarding finite frame theory, we follow the notation in \cite{casazza2013introduction} and refer to it for an explanation of all terms considered here. The basic idea of finite frames is that a redundant system of vectors spanning a vector space can be used in a similar way like a (non-redundant) basis. \neww{In this context, having a redundant basis can be desirable because this means a more robust representation: If information about certain (redundant) components gets lost, e.g. during transmission of data, it is still possible to recover most information. In contrast, a basis representation is, by definition, a minimally viable descriptions and loss of information is catastrophic for recovery.}

In the context of the linear EK, we interpret the vectors 
\begin{equation}\label{eq:finiteframe}
    \varphi^{(i)} := (J-1)^{-1/2} (u^{(i)} - \mu)
\end{equation} as the finite frame vectors, and we can describe any point $p\in \X$ via linear combinations of these vectors. We define the \textit{analysis operator} $T: \X\to\R^J$ and the \textit{synthesis operator} $T^\star: \R^J\to \X$ via
\begin{align}
    T(u) &= ( \langle u, \varphi^{(i)}\rangle )_{i=1}^J\\
    T^\star(\alpha) &= \sum_{i=1}^J \alpha_i \varphi^{(i)}
\end{align}
This gives rise to the \textit{frame operator} $S = T^\star \circ T: \X\to\X$ and the \textit{Grammian operator} $G = T\circ T^\star:\R^J\to \R^J$
\begin{align}
    S(u) &= \sum_{i=1}^J \langle u, \varphi^{(i)}\rangle  \varphi^{(i)}\\
    (G)_{i,j} &= \langle \varphi^{(i)}, \varphi^{(j)}\rangle,
\end{align}
where we identified the linear mapping $G$ with its matrix representation. We note that $S = C_{u,u}$.

\begin{lem}[Reconstruction formula for finite frames, [\cite{casazza2013introduction}]]
    Any point $u\in \X$ can be reconstructed (or spanned) via any of the following four characterisations
\begin{equation}\label{eq:reconstr}
\begin{split}
    u &= \sum_{i=1}^J \langle u, \varphi^{(i)}\rangle S^{-1} \varphi^{(i)}= \sum_{i=1}^J \langle u, S^{-1}\varphi^{(i)}\rangle  \varphi^{(i)}\\
    &= \sum_{i=1}^J \langle S^{-1}u, \varphi^{(i)}\rangle  \varphi^{(i)} = \sum_{i=1}^J \langle u, S^{-1/2}\varphi^{(i)}\rangle  S^{-1/2}\varphi^{(i)}.
\end{split}
\end{equation}
\end{lem}

\subsection{Finite frame bundles for locally weighted ensembles}
When considering locally weighted ensembles, frame vectors become vector bundles $\varphi^{(i)}: \X \to \X$, and the operators become operator fields, e.g., $S: \X\to (\X\to \X)$:
\begin{align}
     \varphi^{(i)}(x) := \sqrt{\kappa^i(x)} (u^{(i)} - \mu^\kappa(x))\\
     S(x)(v) := \sum_{i=1}^J \langle v, \varphi^{(i)}(x)\rangle  \varphi^{(i)}(x)
\end{align}
\begin{lem} The frame operator field is equal to the locally weighted covariance operator field:
    \begin{equation}
        S = C_{u,u}^\kappa
    \end{equation}
    \begin{proof}
        Follows from direct inspection.
    \end{proof}
\end{lem}

\subsection{Local weighting and localisation}
At this point we need to point out that the ideas of local weighting presented in this manuscript are entirely unrelated to the concept of ``localisation'' in the context of Ensemble Kalman (and other) filters. But since this issue is not entirely straightforward, we briefly sketch the demarcation line between these two concepts.

\begin{figure}
    \centering
    \includegraphics[width=\linewidth]{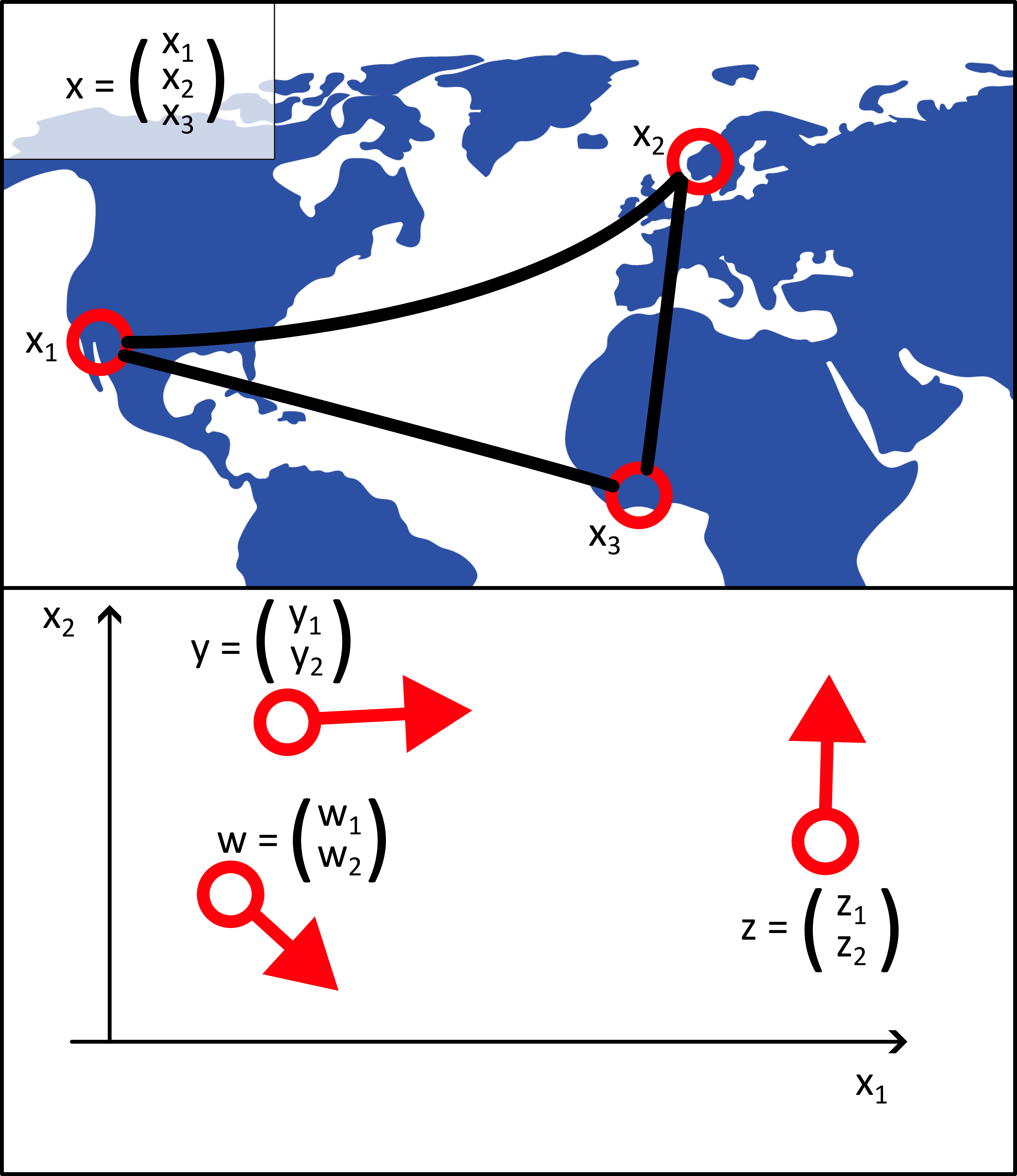}
    \caption{Illustration of Localisation (top) vs local weighting (bottom). Localisation forces decorrelation of components variables $x_i$ in a given dataset where these variables correspond to fixed geographical (or spatial) locations. Local weighting tracks the (changing) spatial configuration of dynamic particles, with ensemble-based effects (such as contraction to a joint mean) being mediated by the particles' respective distance (or a different kernel).}
    \label{fig:illustration_localisation}
\end{figure}

Both traditional localisation and local weighting try to reduce ``effects over long distances'', but the meaning of ``location/distance'' is very different. 
\begin{itemize}
    \item Traditional localisation enforces additional \textbf{decorrelation between parameter entries} corresponding to spatially or geographically distant locations. Its goal is to take into account physical distance $d(i,j)$ between \textit{components} $i,j=1,\ldots,N$, where $\mathbb R^N$ is the parameter space. For instance, let $x = (x_R,x_N,x_T)\in\mathbb R^3$ be the parameter of interest, with $x_R = \text{Temperature in Rome}$, $x_N = \text{Temperature in Naples}$ and $x_T = \text{Temperature in Tokyo}$. Then clearly, due to geographical and meterological domain knowledge we can argue that a-posteriori, $x_R$ and $x_N$ should be more strongly correlated than $x_R$ and $x_T$. Localisation is a tool to enforce decorrelation between $x_R$ and $x_T$.
    \item Local weighting is a tool for creating additional modelling/regression flexibility, by \textbf{tapering off the effect of particles} on each other depending on their (parameter space) distance to each other. This means we penalize distance between \textit{particles} $d(u^{(k)},u^{(l)})$, even if there is no geographical context for them.
    \item Localisation is meaningful only if parameter dimensions are associated with spatial locations.
    \item Localisation is meaningful for problems where the \textit{dimension of parameter space} is at least two (otherwise there are no variables to decorrelate). Local weighting is meaningful if there are at least three \textit{particles} (even in a one-dimensional parameter space).
    \item Localisation takes into account spatial distance $d(i,j)$ between \textit{coordinates} $i$ and $j$ and is computed only once. \new{Local weighting} is based on the distance between particles $d(x^{(k)}, x^{(l)})$ and is iteratively updated since the particles' position evolves in time.
\end{itemize}
\neww{See also figure \ref{fig:illustration_localisation} for a graphical illustration of some of these differences.}

\subsection{Riemannian structure of linear Ensemble Kalman}
 We want to point out that the dynamics of the linear EKI follows an intriguing Riemannian manifold structure. This was remarked in a similar way in the context of affine invariance in \cite{garbuno2020affine}, but it is interesting to make this argument explicitly geometrical. In order to reduce clutter of notation, let $x = u^{(i)}\in \X = \R^d$ be a specific particle following

\begin{equation*}
    \dot x = - C\cdot \nabla \Phi(x),
\end{equation*}
where $C = C^{\kappa}_{u,u}$ (again shortening notation to reduce cluttering).

We can turn this into a geometrically flavored ODE for the components $x^k$ of $x$,
\begin{equation}
    \dot x^k(t) = - M^k_l \delta^{lm} (\nabla \Phi(x(t)))_m
\end{equation}
with $M^k_l = C_{k,l}$ being the components of the $(1,1)$ tensor with the property that $M(\phi, v) = \phi(C\cdot v)$, $\delta^{lm} = 1$ if $l=m$ and $0$ otherwise, and  $(\nabla \Phi(x(t)))_m$ are the components of the co-vector $\nabla \Phi(x(t))$.

Now we define an inverse metric $h^{-1}$ with entries $h^{km} = M^k_l \delta^{lm}$, which means we obtain
\begin{equation}
    \dot x^k(t) = -h^{km} (\nabla \Phi(x(t)))_m = - (\nabla \Phi(x(t)))^k,
\end{equation}
where the last term is what we get if we ``raise the index'' of the covector $\nabla \Phi(x(t))$ using the metric in order to make it a vector. Now we fuse time and space by defining the particle dynamics as a curve in space-time:

We augment the particle trajectory by a time component, $X = (t,x) \in \R^{1+d}$: Defining $\Psi(X) = \Phi(x) - t$ and 
\begin{equation*}
    K(t) = \begin{pmatrix}
        1 & 0 \\ 0 & C
    \end{pmatrix}
\end{equation*} 
(noting that $C$ depends on $t$ via the time-dependent ensemble/mean-field) turns our ODE into
\begin{equation}
    \dot X^k(t) = - K(t)^k_l \delta^{lm} \nabla \Psi(X(t))_m. 
\end{equation}

By defining an inverse metric $g^{km} = K(t)^k_l \delta^{lm}$, and equivalently a metric 
\begin{equation}
     g_{km} = \begin{pmatrix}
        1 & 0 \\ 0 & C(t)^{-1}
    \end{pmatrix}_{k,m},
\end{equation}
this means we have
\begin{equation}
    \dot X^k(t) = -g^{km} \nabla \Psi(x(t))_m = - \nabla \Psi(x(t))^k
\end{equation}

A striking property of this evolution equation is the closed-loop property of the coupled system of ensemble and curvature. Adapting J. A. Wheeler's aphorism about general relativity: \textit{The particles' mean field $\nu$ tells spacetime how to curve. Spacetime tells particles how to move:}

Or, stated visually in terms of an ensemble (which has the same property as the mean-field case):
\begin{equation}
\{u^{(i)}\}_{i=1}^J\tikzmark{b} \stackrel{ C_{u,u}^\kappa}{\mapsto} g_{km} \mapsto \tikzmark{a}\{\dot u^{(i)}\}_{i=1}^J
\begin{tikzpicture}[overlay,remember picture,out=225,in=315,distance=0.5cm]
    \draw[->,shorten >=3pt,shorten <=3pt] (a.center) to (b.center);
  \end{tikzpicture}
\end{equation}

Another way of thinking about this is that the EKI (both in its linear and locally weighted formulation) performs a natural gradient flow on the manifold generated by the empirical covariance.

\section{Ensemble-based linear approximation}\label{sec:linearapprox}
In this section we flesh out a notion of ensemble-based linearisation based on local weighting. This will allow us to interpret the role of the locally weighted versions of moments in \eqref{eq:locally_weighted_moments} in the definition of locally weighted EK \eqref{eq:lwEK}. We provide a brief motivation for local linear approximation by thinking about the action of $A$ in terms of the reconstruction formula and by approximating derivatives of $A$ via ensemble evaluations $A(u^{(i)})$. There is a connection of the ideas presented here to locally weighted polynomial regression (LOESS/LOWESS), see \cite{cleveland1979robust}, but the point of view and scope vary between these approaches. In some sense, LOESS implements $\mu_A^\kappa(x)$ as a function of $x$, while we are interested in a linear approximation centered at (fixed values of) $x$, correct up to first or second order of approximation. Simplex gradients (see \cite{bortz1998simplex}, \cite{conn2009introduction}, \cite{coope2019efficient}, \cite{custodio2008using}) are a related idea from gradient-free optimisation, but don't contain weights based on distance. 

Another point of view on what follows is that ensemble-based gradient approximation acts as a statistical analogue to symmetrised difference quotients used for approximating derivatives.

We start by first building approximate derivatives based on ensembles or their generated mean fields and leverage this for function approximation in section \ref{sec:func_approx}.
\subsection{Ensemble-based approximate derivatives}
In this section we demonstrate how pointwise evaluations $\{A(u^{(i)})\}_{i=1}^J$ can be used to construct an approximation to the derivatives $DA(x)[f]$ and $D^2A(x)[f,g]$ of $A$. Prior work (in \cite{schillings2023ensemble}) investigated this problem via optimisation of a least-squares functional, but here we show that this can be done in a much more direct way by constructing weighted empirical covariance operators.

Let $A: \X\to \Y$ be a linear map and $\{ u^{(i)}-\mu\}_{i=1}^J$ be a finite frame for $\X$. Then, due to the reconstruction formula \eqref{eq:reconstr},
\begin{align*}
    Av = \sum_{i}A( u^{(i)}-\mu)\langle  u^{(i)}-\mu, S^{-1} v\rangle = C_{A,u}C_{u,u}^{-1}v,
\end{align*}
which means, in the linear case, $A = C_{A,u}C_{u,u}^{-1} = DA(v)$, where $DA$ is the Jacobian of the mapping $A$. Of course this is just a needlessly complicated way of writing $Av$, but these ideas suggest the following definition of a frame-based derivatives of nonlinear mapping $A$, allowing for incorporation of local weights $\kappa^{(i)}$.
\begin{definition}\label{def:frame_dep_gradient}
    Let $\kappa:=\{\kappa^{(i)}: \X\to \R\}_{i=1}^J$ be a set of weight fields summing up to $1$. Let $A: \X\to \Y$ be a (possibly nonlinear) function. Then the ensemble-based derivative of $A$ is defined as the mapping $D_\kappa  A: \X \to \mathcal L(\X,\Y)$
    \begin{equation}
    \begin{split}
        D_\kappa  A(x)[f] &:= \sum_{i=1}^J \kappa^{(i)}(x) (A(u^{(i)}) - \mu_A^\kappa(x)) \langle u^{(i)} - \mu^\kappa(x), C_{u,u}^\kappa(x)^{-1}[f]\rangle \\
         &= (C_{A,u}^\kappa(x)\circ C_{u,u}^\kappa(x)^{-1})[f]
    \end{split}
    \end{equation}
\end{definition}
\begin{lem}
    An equivalent definition is 
    \begin{equation}
    \begin{split}
         D_\kappa  A(x)[f] &:= \sum_{i=1}^J \kappa^{(i)}(x) (A(u^{(i)}) - A(\mu^\kappa(x))) \langle u^{(i)} - \mu^\kappa(x), C_{u,u}^\kappa(x)^{-1}[f]\rangle
    \end{split}
    \end{equation}
    (Note that we replaced $\mu_A^\kappa(x)$ by $A(\mu^\kappa(x))$). 
    \begin{proof}This is true because 
    \begin{align*}
        &\sum_{i=1}^J \kappa^{(i)}(x) (\mu_A^\kappa(x) - A(\mu^\kappa(x))) \langle u^{(i)} - \mu^\kappa(x), C_{u,u}^\kappa(x)^{-1}[f]\rangle\\
        &= (\mu_A^\kappa(x) - A(\mu^\kappa(x))) \left\langle \sum_{i=1}^J \kappa^{(i)}(x)  ( u^{(i)} - \mu^\kappa(x)), C_{u,u}^\kappa(x)^{-1}[f]\right\rangle = 0
    \end{align*}
    by \eqref{eq:fundamental}.\end{proof}
\end{lem}

To save some notation, we will from now on always understand $C=C_{u,u}^\kappa(x)$, where the choice of $x$ will be clear from context.

It is interesting to further define a notion of second-degree frame-dependent derivative as follows.

\begin{definition}
    In the setting of Definition \ref{def:frame_dep_gradient}, we set
\begin{align*}
   D^2_\kappa A(x)[f,g] &= \sum_{ij}\kappa^i(x)\kappa^j(x)  \langle u^{(i)} - \mu^\kappa(x), C^{-1}f\rangle \langle u^{(j)} - \mu^\kappa(x), C^{-1}g\rangle \\
    &\cdot (D_\kappa  A(u^{(i)}) - D_\kappa  A(\mu^\kappa(x)))[u^{(j)} - \mu^\kappa(x)]
\end{align*}
\end{definition}

This expression seems to be asymmetric in the summation indices $i$ and $j$ so some justification is needed:

\begin{lem}\label{lem:symmetrised}
Firstly, $D^2_\kappa A(x)[f,g] = D_\kappa (x\mapsto D_\kappa  A(x)[f])[g]$ (which means $D_\kappa ^2$ deserves the symbol it has). 

Secondly, if we define the symmetrised variant
    \begin{align*}
    M_0 &:= \sum_{ij}\kappa^i(x)\kappa^j(x)  \langle u^{(i)} - \mu^\kappa(x), C^{-1}f\rangle \langle u^{(j)} - \mu^\kappa(x), C^{-1}g\rangle \\
    &\cdot \left\{ \frac{(D_\kappa  A(u^{(j)}) - D_\kappa  A(\mu^\kappa(x)))[u^{(i)} - \mu^\kappa(x)]}{2} \right.+\\ &+\left.\frac{(D_\kappa  A(u^{(i)}) - D_\kappa  A(\mu^\kappa(x)))[u^{(j)} - \mu^\kappa(x)]}{2}\right\}
\end{align*}
 and the antisymmetrised variant
 \begin{align*}
    M_2 &:= \sum_{ij}\kappa^i(x)\kappa^j(x)  \langle u^{(i)} - \mu^\kappa(x), C^{-1}f\rangle \langle u^{(j)} - \mu^\kappa(x), C^{-1}g\rangle \\
    &\cdot (D_\kappa  A(u^{(j)}) - D_\kappa  A(\mu^\kappa(x)))[u^{(i)} - \mu^\kappa(x)]
\end{align*}
Then by definition, $M_1 := D^2_\kappa A(x)[f,g] = 2\cdot M_0 - M_2$, but more importantly, 
\[D^2_\kappa A(x)[f,g] = M_0 = M_2,\] 
which means we can use any of the variants derived as alternative definitions of $D^2_\kappa A(x)[f,g]$.
\begin{proof}
    For brevity, we sometimes write $\varphi^{(j)}(x) = u^{(j)} - \mu^\kappa(x)$ for $j=1,\ldots, J$. 

Regarding the first statement, using the reconstruction formula 
    \begin{align*}
        &\sum_j \kappa^{(j)}(x) C^{-1}(u^{(j)}-\mu^\kappa(x)) \left\langle D_\kappa  A(u^{(i)})  - D_\kappa  A(\mu^\kappa(x)), u^{(j)}-\mu^\kappa(x) \right\rangle\\
        &= D_\kappa  A(u^{(i)}) - D_\kappa  A(\mu^\kappa(x)),
    \end{align*} 
    we get
    \begin{align*}
        &D_\kappa (x\mapsto D_\kappa  A(x)[f])[g]\\
        &= \sum_i \kappa^{(i)}(x)  (D_\kappa  A(u^{(i)})[f] - D_\kappa  A(\mu^\kappa(x))[f]) \cdot \langle \varphi^{(i)}(x)), C^{-1}g\rangle \\
        &=\sum_{i} \kappa^{(i)}(x) \left\langle f, \sum_j \kappa^{(j)}(x) C^{-1}\varphi^{(j)}(x) \langle D_\kappa  A(u^{(i)})  - D_\kappa  A(\mu^\kappa(x)),\varphi^{(j)}(x)\rangle\right\rangle\\
        &\qquad\cdot \langle \varphi^{(i)}(x)), C^{-1}g\rangle \\
        &=\sum_{i,j} \kappa^{(i)}(x)\kappa^{(j)}(x) \left( D_\kappa  A(u^{(i)})  -D_\kappa  A(\mu^\kappa(x))\right)[\varphi^{(j)}(x)]\\
        &\quad \cdot \langle \varphi^{(j)}(x)), C^{-1}f\rangle \cdot \langle \varphi^{(i)}(x)), C^{-1}g\rangle =D_\kappa ^2 A(x)[f,g].
    \end{align*}

Regarding the alternative representation claim: We start with $D_\kappa^2 A(x)[f,g]$:
    \begin{align*}
        M_1&=\sum_{ij}\kappa^i(x)\kappa^j(x)  \langle \varphi^{(i)}(x), C^{-1}f\rangle \langle \varphi^{(j)}(x), C^{-1}g\rangle \\
    &\qquad\cdot (D_\kappa  A(u^{(i)}) - D_\kappa  A(\mu^\kappa(x)))[\varphi^{(j)}(x)]\\
    &=\sum_{ij}\kappa^i(x)\kappa^j(x)  \langle \varphi^{(i)}(x), C^{-1}f\rangle \langle \varphi^{(j)}(x), C^{-1}g\rangle \\
    &\qquad\cdot (D_\kappa  A(u^{(i)}) - D_\kappa  A(u^{(j)}))[\varphi^{(j)}(x)]\\
    &+\sum_{ij}\kappa^i(x)\kappa^j(x)  \langle \varphi^{(i)}(x), C^{-1}f\rangle \langle \varphi^{(j)}(x), C^{-1}g\rangle \\
    &\qquad\cdot (D_\kappa  A(u^{(j)}) - D_\kappa  A(\mu^\kappa(x)))[u^{(j)} - u^{(i)}]\\
    &+\sum_{ij}\kappa^i(x)\kappa^j(x)  \langle \varphi^{(i)}(x), C^{-1}f\rangle \langle \varphi^{(j)}(x), C^{-1}g\rangle \\
    &\qquad\cdot (D_\kappa  A(u^{(j)}) - D_\kappa  A(\mu^\kappa(x)))[\varphi^{(i)}(x)]\\
    &=: I_1 + I_2 + M_2
    \end{align*}
    We take a closer look at $I_1$.
    \begin{align*}
        I_1 &= \sum_{ij}\kappa^i(x)\kappa^j(x)  \langle \varphi^{(i)}(x), C^{-1}f\rangle \langle \varphi^{(j)}(x), C^{-1}g\rangle \\
    &\qquad\cdot (D_\kappa  A(u^{(i)}) - D_\kappa  A(u^{(j)}))[\varphi^{(j)}(x)]\\
    &=\sum_j \kappa^j(x) \langle \varphi^{(j)}(x), C^{-1}g\rangle \\
    &\qquad\left[\sum_i \kappa^i(x)  (D_\kappa  A(u^{(i)}) - D_\kappa  A(u^{(j)})) \langle \varphi^{(i)}(x), C^{-1}f\rangle  \right](\varphi^{(j)}(x))\\
    &=\sum_j \kappa^j(x) \langle \varphi^{(j)}(x), C^{-1}g\rangle \\
    &\qquad\underbrace{\left[\sum_i \kappa^i(x)  (D_\kappa  A(u^{(i)}) - D_\kappa  A(\mu^\kappa(x)))[\varphi^{(j)}(x)] \langle \varphi^{(i)}(x), C^{-1}f\rangle \right]}_{D_\kappa D_\kappa  A(x)[\varphi^{(j)}(x)][f]}\\
    &=D_\kappa D_\kappa  A(x)\underbrace{\left[\sum_j \kappa^j(x) \langle \varphi^{(j)}(x), C^{-1}g\rangle (\varphi^{(j)}(x))\right]}_{g}[f] \\
    &= D_\kappa D_\kappa  A(x)[f,g]
    \end{align*}
    Term $I_2$ is
    \begin{align*}
        I_2 &= \sum_{ij}\kappa^i(x)\kappa^j(x)  \langle \varphi^{(i)}(x), C^{-1}f\rangle \langle \varphi^{(j)}(x), C^{-1}g\rangle \\
    &\qquad\cdot (D_\kappa  A(u^{(j)}) - D_\kappa  A(\mu^\kappa(x)))[u^{(j)} - u^{(i)}]\\
    &= \sum_j \kappa^j(x) \langle \varphi^{(j)}(x), C^{-1}g\rangle (D_\kappa  A(u^{(j)}) - D_\kappa  A(\mu^\kappa(x)))\\
    &\qquad \underbrace{\left[\sum_i \kappa^i(x) \langle \varphi^{(i)}(x), C^{-1}f\rangle (u^{(j)} - u^{(i)}) \right]}_{=0-f}\\
    &= \sum_j \kappa^j(x) \langle \varphi^{(j)}(x), C^{-1}g\rangle (D_\kappa  A(u^{(j)}) - D_\kappa  A(\mu^\kappa(x)))[-f]\\
    &= -D_\kappa D_\kappa  A(x)[f,g].
    \end{align*}
    This shows that $I_1 + I_2 = 0$, i.e. $M_1 = M_2$. 
\end{proof}
\end{lem}

\subsection{Mean-field-based approximate derivatives}
All notions considered so far generalise in a natural way to the setting $J\to\infty$, i.e. instead of an ensemble $\mathcal E = \{u^{(i)}\}_i$ we consider a continuous (probability) measure $\nu$: A choice of $\nu$-integrable kernel $k$ gives rise to a local weighting function
\begin{equation}\label{eq:weighting_func}
    \kappa(x, z) := \frac{k(x,z)}{\int_\X k(x, w) \d\nu(w)},
\end{equation}
local mean fields
\begin{align}
    \underline{\mu}^\kappa(x) &:= \int_\X z \kappa(x,z) \d\nu(z)\\
    \underline\mu_A^\kappa(x) &:= \int_\X A(z) \kappa(x,z) \d\nu(z)
\end{align}
as well as covariance forms
\begin{align}
    \underline C_{u,u}^\kappa(x) &:= \int_\X \kappa(x,z) (z - \underline\mu^\kappa(x))\cdot (z-\underline\mu^\kappa(x))^\top \d\nu(z)\\
    \underline C_{u,A}^\kappa(x) &:= \int_\X \kappa(x,z) (z - \underline\mu^\kappa(x))\cdot (A(z)-\underline\mu_A^\kappa(x))^\top \d\nu(z).
\end{align}

The mean-field-based derivative of $A$ is defined as the mapping
   \begin{equation}
    \begin{split}
        \underline D_\kappa  A(x)[f] &:= \int_\X \kappa(x,z) (A(z) - \underline\mu_A^\kappa(x)) \langle z - \underline\mu^\kappa(x), \underline C_{u,u}^\kappa(x)^{-1}[f]\rangle  \d\nu(z)
    \end{split}
    \end{equation}
with a similar expression for $\underline D_\kappa ^2 A(x)[f,g]$.

\begin{rem}
    While the ensemble-form and the mean-field form of $D_\kappa  A$ are certainly different, it will usually be clear from context which version we mean.
\end{rem}

\subsection{Ensemble-based function approximation}\label{sec:func_approx}
Having defined approximate ensemble-based and mean-field-based derivatives $D_\kappa  A(x)[f]$ and $D_\kappa  A(x)[f]$ we can define the following notions of local approximation (to first and second order) of the mapping $A$, from the point of view of $x\in \X$, utilising the ensemble $\{u^{(i)}\}_{i=1}^J$ and kernel weights $\kappa^i$ (or the mean field $\nu$). \new{An important design choice for approximating nonlinear maps with an ensemble is the order in which \textit{empirical averages} and \textit{nonlinear map application} are performed: Both $\mu_A^\kappa(x)$ (locally weighted empirical average of values of $A(x)$), and $A(\mu^\kappa(x))$ (evaluation of map $A$ in locally weighted empirical average) are sensible choices leading to slightly different approximations.}
\begin{align*}
     A_x(\xi) &:= \mu_A^\kappa(x) + D_\kappa  A(x) \cdot (\xi - \mu^\kappa(x))\\
    \widetilde A_x(\xi) &:= A(\mu^\kappa(x)) + D_\kappa  A(x) \cdot (\xi - \mu^\kappa(x))\\
     A_x^{(2)}(\xi) &:= \mu_A^\kappa(x) + D_\kappa  A(x) \cdot (\xi - \mu^\kappa(x)) + \frac{1}{2}D_\kappa ^2A(x)[\xi - \mu^\kappa(x),\xi - \mu^\kappa(x)]\\
     \widetilde A_x^{(2)}(\xi) &:=  A(\mu^\kappa(x)) + D_\kappa  A(x) \cdot (\xi - \mu^\kappa(x)) + \frac{1}{2}D_\kappa ^2A(x)[\xi - \mu^\kappa(x),\xi - \mu^\kappa(x)],
\end{align*}
where the mean-field versions $\underline A_x, \underline {\widetilde A}_x,\underline A_x^{(2)}, \underline {\widetilde A}_x^{(2)}$ of these expressions are straightforward by replacing $\mu^\kappa \mapsto \underline \mu^\kappa$ etc.

We prove the next Lemma\new{, which illustrates the most important properties of (and distinctions between) these variants, constraining ourselves} for a specific Gaussian-type kernel to simplify calculations. We do note that these statements can be generalised to a much broader class of kernels satisfying sufficient integrability.

\begin{lem}[Properties of the ensemble approximation]
We choose a Gaussian-type kernel with bandwidth $r$, i.e. $k(x,z) = {(2\pi)^{-d/2} r^{-d} }\exp(-\|x-z\|^2/(2r))$. Then the linear approximations $A_x$ and $\widetilde A_x$ have the following properties.
\begin{enumerate}
    \item $\widetilde A_x(\mu^\kappa(x)) = A(\mu^\kappa(x))$ 
    \item For fixed $x\in \mathcal X$, the mapping $ A_x$ is the minimiser of the least squares functional $V$, where
        \begin{equation}
            V(L) = \frac{1}{2}\sum_i \kappa^i(x) \|A(u^{(i)}) - L(u^{(i)})\|^2
        \end{equation}
        among all affine-linear regression maps $\mathcal{A} = \{L: \mathcal X\to \mathcal Y,~L(\xi) = M_x(\xi-\mu^\kappa(x)) + b_x\}$., i.e.
    \begin{equation}
        \argmin_{L \in \mathcal{A}} V(L) = A_x
    \end{equation}
    In particular, the optimal choice corresponds to $M_x = D_\kappa  A(x)$ and $b_x = \mu_A^\kappa(x)$,
    \item The estimators $\widetilde A_x$ and $A_x$ satisfy a bias-variance trade-off: $V(\widetilde A_x) = V(A_x) + \|A(\mu^\kappa(x)) - \mu_A^\kappa(x)\|^2$

\end{enumerate}
\begin{proof}[Sketch of proof]
The key technical arguments are standard applications of Taylor's theorem and the Laplace approximation (\cite{wong2001asymptotic}), so we just give a sketch of the proof. The Laplace approximation allows to show that \begin{equation}\int h(z)\kappa(x,z)\d\nu(z) \sim h(x)\end{equation} for $r\to 0$. The order of asymptotics can be found in \cite{schillings2020convergence}, but this point is of secondary importance here.

    1. \new{$\widetilde A_x(\mu^\kappa(x)) = A(\mu^\kappa(x)) + D_\kappa  A(x) \cdot (\mu^\kappa(x) - \mu^\kappa(x)) = A(\mu^\kappa(x)) + 0$.} 
    
    2.: We consider the functional 
    \[V(b_x, M_x) = \frac{1}{2}\sum_i \kappa^i(x) \|A(u^{(i)}) - M_x(u^{(i)} - \mu^\kappa(x)) - b_x\|^2.\]
    Standard differentiation rules (see, e.g., \cite{petersen2008matrix}) yield
    \begin{align*}
        \partial_{b_x}V &= -\sum_i \kappa^{(i)}(x) (A(u^i) - b_x - M_x(u^{(i)}-\mu^\kappa(x))) \\
        &= -\sum_i \kappa^{(i)}(x) (A(u^i) - b_x)\\
        \partial_{M_x}V &= -\sum_i \kappa^{(i)}(x) (A(u^i) - b_x - M_x(u^{(i)}-\mu^\kappa(x)))(u^{(i)} - \mu^\kappa(x))^\top \\
        &= -\sum_i \kappa^{(i)}(x) (A(u^i) - b_x)(u^{(i)} - \mu^\kappa(x))^\top  \\
        &\qquad+ M_x\sum_i \kappa^{(i)}(x) (u^{(i)}-\mu^\kappa(x)))(u^{(i)} - \mu^\kappa(x))^\top 
    \end{align*}
    Setting both partial derivatives to $0$ yields first $b_x = \mu_A^\kappa(x)$ and then $M_x = D_\kappa  A(x)$, as claimed.

    3.: Follows from straightforward calculation.
\end{proof}
\end{lem} 

\begin{lem}[Properties of the mean field approximation]
    We choose a Gaussian-type kernel with bandwidth $r$, i.e. $k(x,z) = {(2\pi)^{-d/2} r^{-d} }\exp(-\|x-z\|^2/(2r))$. We assume that the mean field $\nu$ is absolutely continuous with respect to Lebesgue measure. Then the mean-field linear approximations $\underline A_x$ and $\underline{\widetilde A}_x$ have the following properties.
    \begin{enumerate}
        \item When the bandwidth $r\searrow 0$ then $\widetilde A_x(\xi)$ is the first-order Taylor approximation at $(x,A(x))$.
        \item When the bandwidth  $r\to \infty$, 
        \begin{itemize}
            \item $D_\kappa  A_x = C_{A,u}C_{u,u}^{-1}$ (without any weighting),
            \item $A_x$ is the usual ordinary least-squares affine-linear regression function, and
            \item $A_x^{(2)} = 0$.
        \end{itemize} 
    \end{enumerate}
\begin{proof}
    1. can be proven with the Laplace approximation:
    \begin{align*}
        \lim_{r\searrow 0}\mu^{\kappa_r}(x) &= \lim_{r\searrow 0}\int z \kappa_r(x,z) \d\nu(z) = x \\
        \lim_{r\searrow 0} \mu^{\kappa_r}_A(x) &= \lim_{r\searrow 0} \int A(z) \kappa_r(x,z) \d\nu(z)\\
        &= \int \left[A(x) + DA(x)(z-x) + \mathcal O(\|x-z\|^2) \right]\kappa_r(x,z)\d\nu(z)\\
        &=A(x) + DA(x)(\lim_{r\searrow 0} \mu^{\kappa_r}(x)-x) + \int  \mathcal O(\|x-z\|^2)\kappa_r(x,z)\d\nu(z)\\
        &= A(x)
    \end{align*}
    Also, sketching the argument using Taylor's theorem,
    \begin{align*}
        \lim_{r\searrow 0}D_{\kappa_r}  A(x) &= \int_\X \kappa_r(x,z) (A(z) - \mu^{\kappa_r}_{A}(x)) (z - \mu^{\kappa_r}(x))^\top   \d\nu(z) C_{u,u}^{\kappa_r}(x)^{-1}\\
        &\simeq \int_\X \kappa_r(x,z) (A(z) - A(x)) (z - x)^\top   \d\nu(z) C_{u,u}^{\kappa_r}(x)^{-1} \\
        &\simeq \int_\X \kappa_r(x,z) \left[DA(x)\cdot(z-x) + \mathcal O(\|z-x\|^2)\right] (z - x)^\top   \d\nu(z) C_{u,u}^{\kappa_r}(x)^{-1} \\
        &\simeq DA(x) C_{u,u}^{\kappa_r}(x) \cdot C_{u,u}^{\kappa_r}(x)^{-1} = DA(x)
    \end{align*}
    2. follows from straightforward calculation.
\end{proof}
\end{lem}

\begin{rem}
    Note that $A_x$ and $\widetilde A_x$ only differ by the constant shift $\mu_A^\kappa(x) - A(\mu^\kappa(x))$. $A_x$ is the least-square minimising affine-linear approximation, while $\widetilde A_x$ preserves the pointwise property $\widetilde A_x(\mu^\kappa(x)) = A(\mu^\kappa(x))$. This also holds in an equivalent way for the mean-field versions of these objects.
\end{rem}
\neww{
\begin{rem}
    While the analysis of how well a given ensemble can approximate a map $A$ via $A_x$ and $A_x^{(2)}$ is interesting and important, sometimes ``bad approximation'' can be beneficial: In the discontinuous example in section \ref{sec:discontinuous}, we see that the exact gradient does not provide any useful information for finding directions for minimisation, but the (inexact) ensemble-based approximation supplies a ``smoothed-out idea'' of this objective function which can be leveraged for finding minima of the function.
\end{rem}
}
\subsection{Examples}
\begin{figure}
    \centering
    \includegraphics[width=\linewidth]{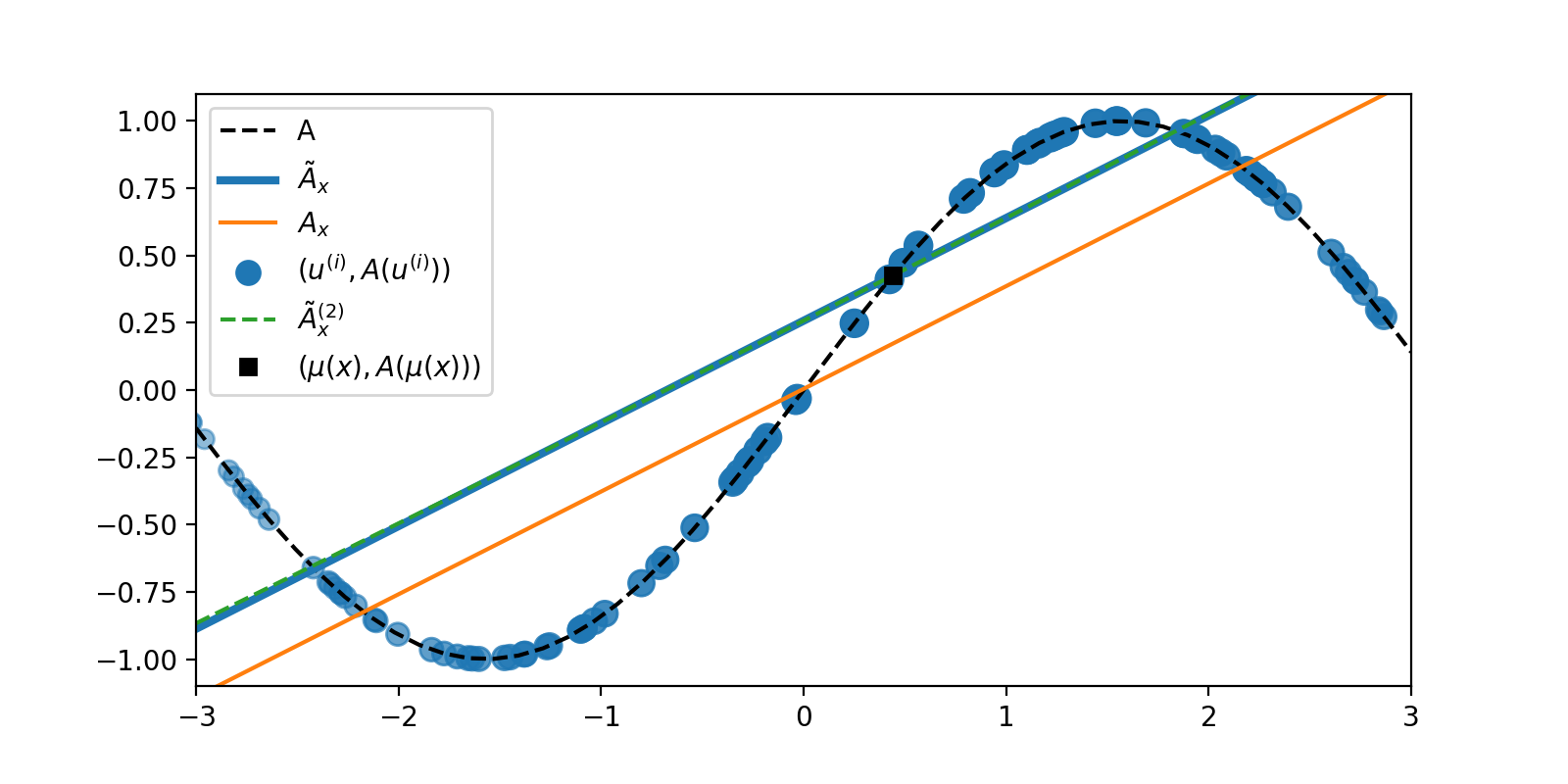}\\
    \includegraphics[width=\linewidth]{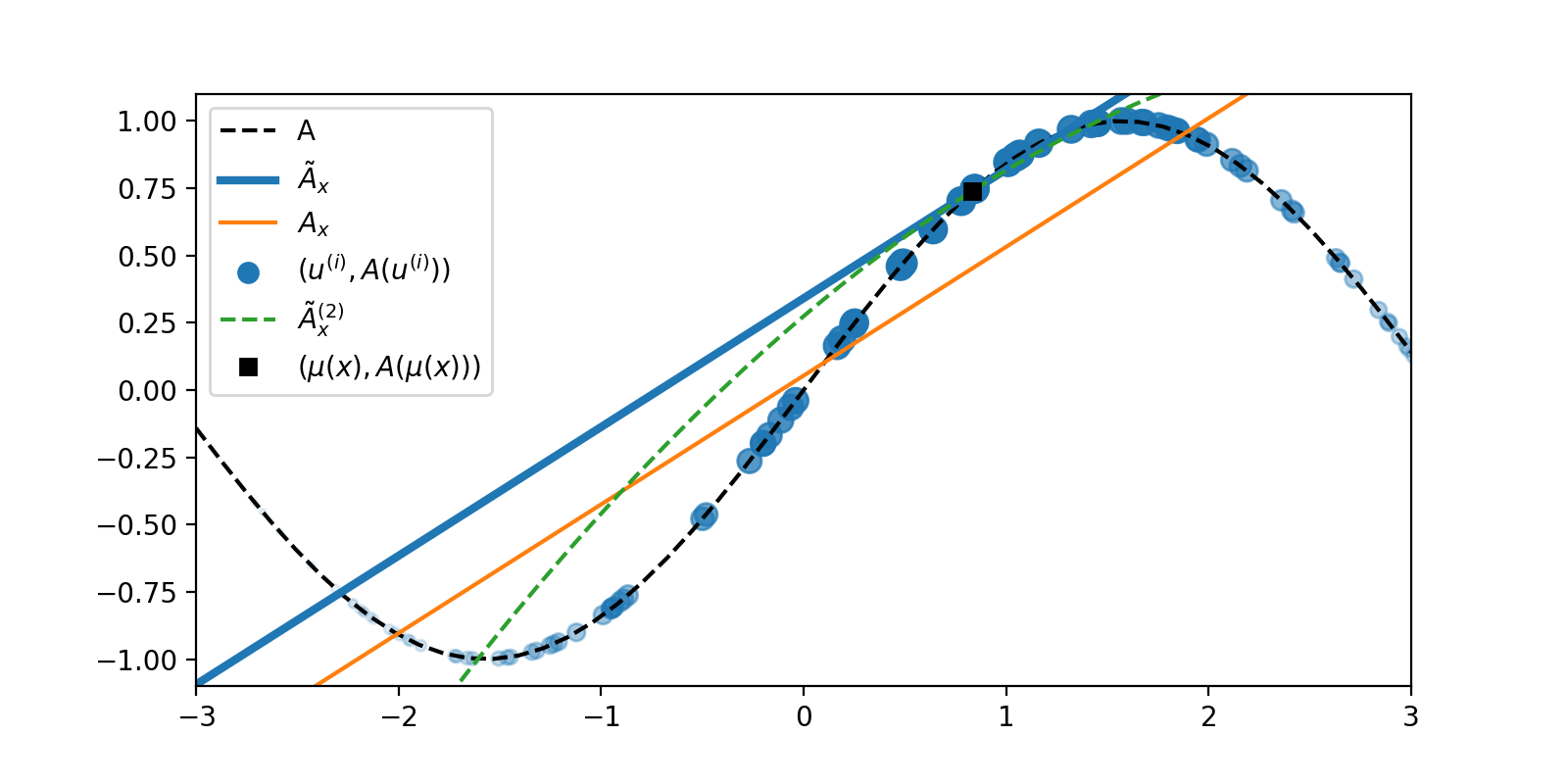}\\
    \includegraphics[width=\linewidth]{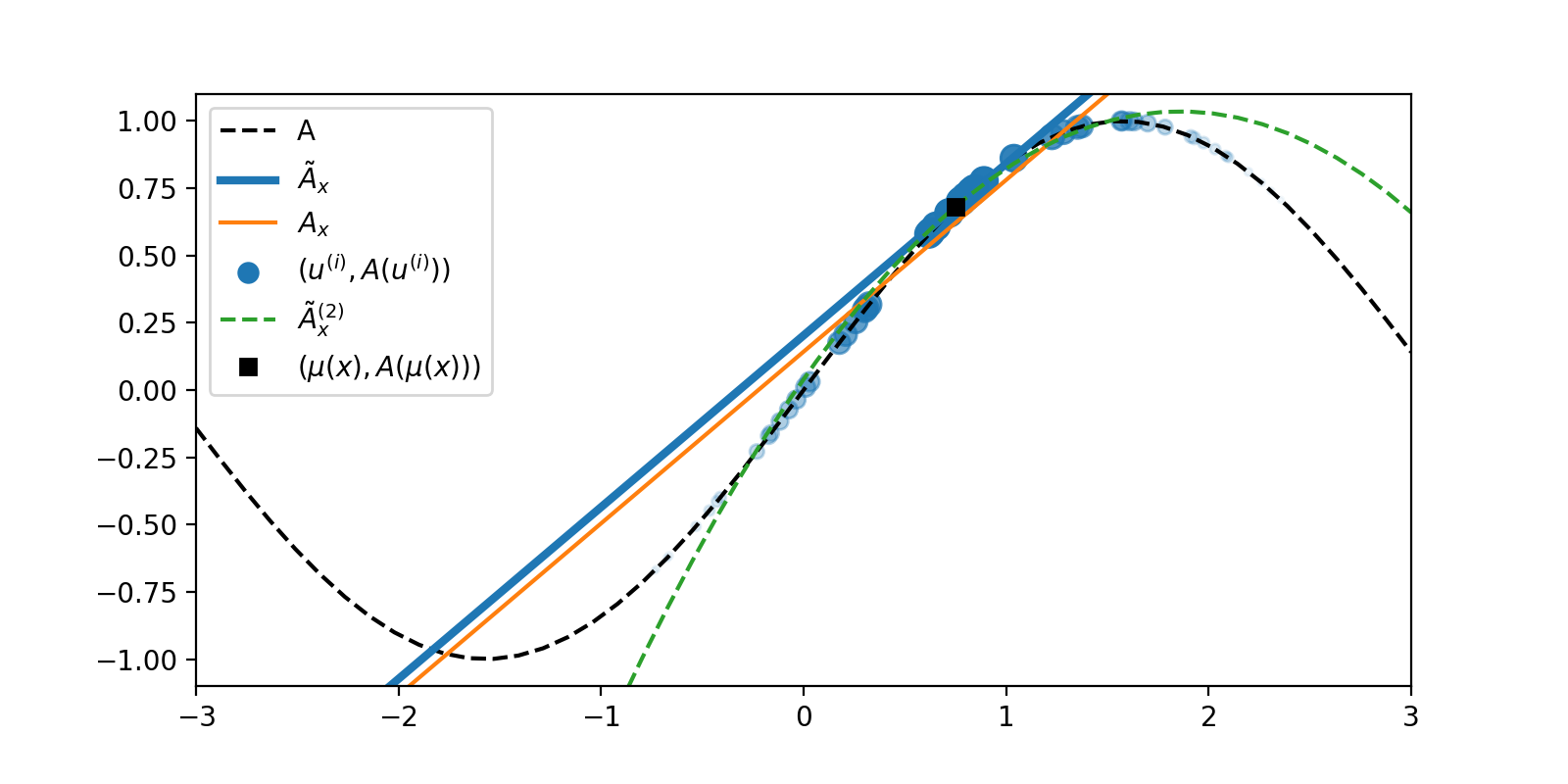}%
    \caption{Local frame-based approximation with (Gaussian) kernel bandwidths $r\in\{5,1,0.2\}$ (top to bottom) for the sine function, with an ensemble of size $J=100$ uniformly sampled on the interval $[-3,3]$. }
    \label{fig:approx_sin}
\end{figure}

\new{We give a few examples to demonstrate the performance of $A_x$ and $A_x^{(2)}$ in approximating a given map $A$.}
\paragraph{Sine function.} As a first example, we choose $\X = \Y = \R$ and $A(x) = \sin(x)$. We generate an ensemble with $J=100$ particles uniformly on the interval $[-2,4]$ and pick a reference point $x=\pi/4$. Figure \ref{fig:approx_sin} shows the result of three different bandwidths, $r\in\{5,1,0.2\}$ for a Gaussian kernel $k(x,y) = \exp(-(x-y)^2/(2r^2))$. \new{The three different plots show three different bandwidth, and their impact on ``locality'' of the approximation. We see how the largest bandwidth, $r=5.0$ (top row of figure \ref{fig:approx_sin}), leads to a global approximation similar to an affine-linear regression fit. The smallest bandwidth $r=0.2$ (bottom row) leads to a Taylor-esque local approximation around the reference point $x=\pi/4$. All three figures show the true function $A$ to be approximated, the two first-order approximation variants $A_x$ and $\tilde A_x$, the second order approximation $\tilde A_x^{(2)}$, and the position of the locally weighted empirical mean $\mu(x)$. Note that the ensemble is the same in every single case (local weighting is visually conveyed by the fact that ensemble points close to the standard deviation of the local weighting kernel are printed with a lower opacity value (and vanish completely beyond that)}

\begin{figure}
    \centering
    \includegraphics[width=\linewidth]{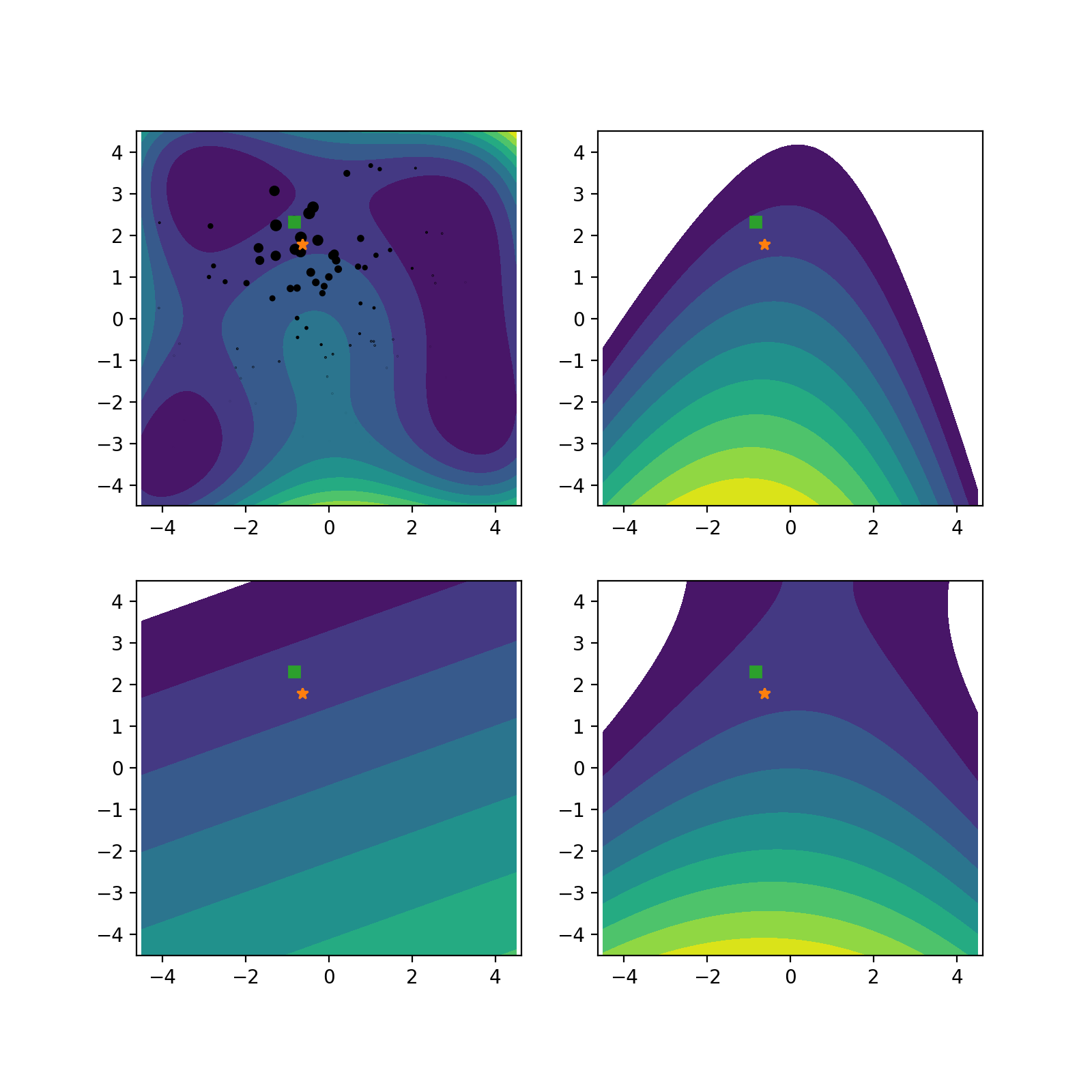}
    \caption{Local frame-based approximation with (Gaussian) kernel bandwidth $r=1$. Ensemble size $J=100$. \textit{Top left:} Contour plot of the Himmelblau function. Green square marker: Reference point $x$. Black markers: Ensemble members (with marker size proportional to weight $\kappa^{(i)}(x)$). Orange star: locally weighted ensemble mean $\mu^\kappa(x)$. \textit{Top right:} Contour plot of (exact) second-order Taylor approximation centered at $\mu^\kappa(x)$. \textit{Bottom left: }Approximate first order approximation $A_x$. \textit{Bottom right: }Approximate second-order approximation $A_x^{(2)}$.}
    \label{fig:approx_himmelblaur1}
\end{figure}

\begin{figure}
    \centering
    \includegraphics[width=\linewidth]{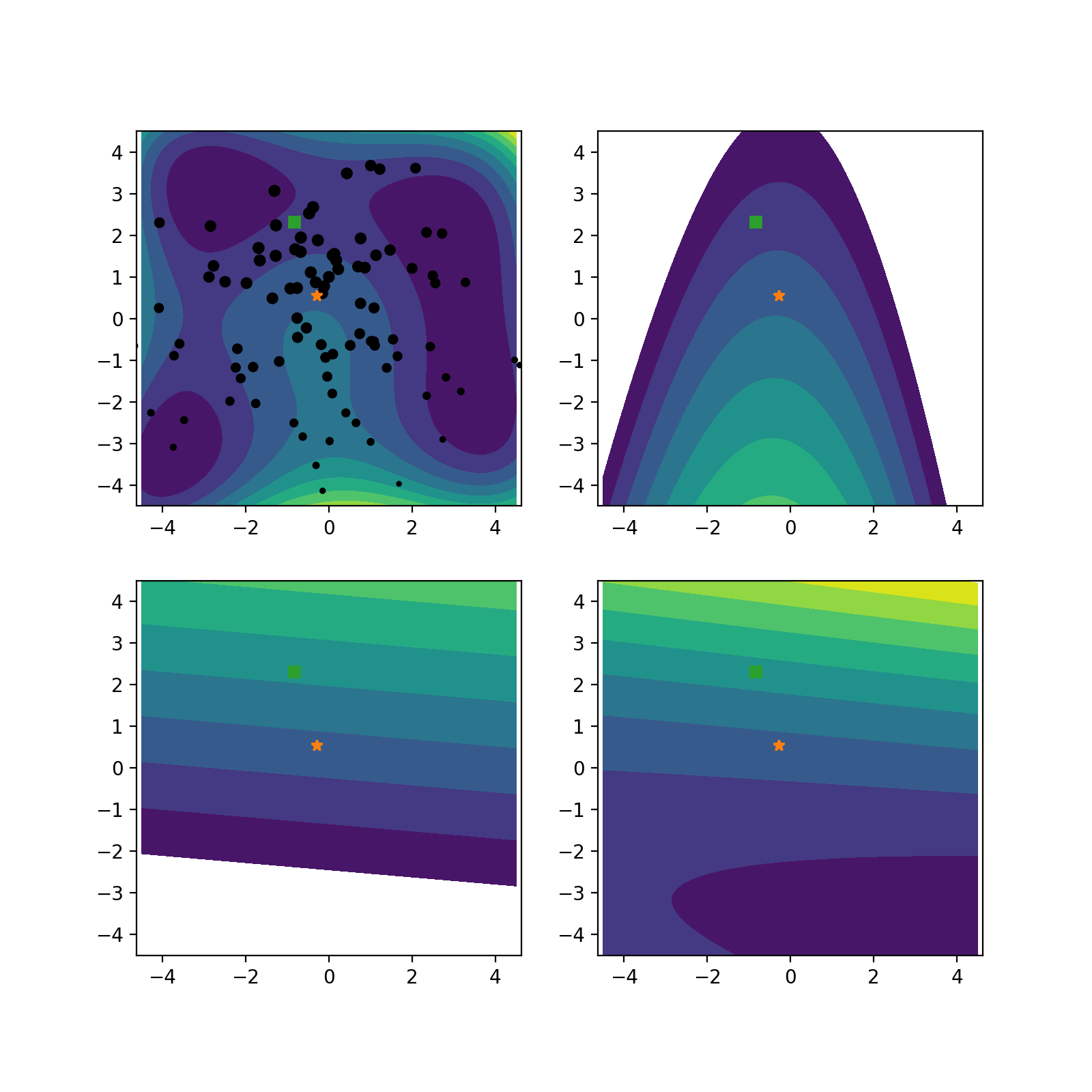}
    \caption{Local frame-based approximation with (Gaussian) kernel bandwidth $r=5$. Ensemble size $J=100$. \textit{Top left:} Contour plot of the Himmelblau function. Green square marker: Reference point $x$. Black markers: Ensemble members (with marker size proportional to weight $\kappa^{(i)}(x)$), identical ensemble as in figure \ref{fig:approx_himmelblaur1}. Orange star: locally weighted ensemble mean $\mu^\kappa(x)$. \textit{Top right:} Contour plot of (exact) second-order Taylor approximation centered at $\mu^\kappa(x)$. Note that this is different from the top right plot in figure \ref{fig:approx_himmelblaur1} because $\mu^\kappa(x)$ depends on the kernel and its bandwidth. \textit{Bottom left: }Approximate first order approximation $A_x$. \textit{Bottom right: }Approximate second-order approximation $A_x^{(2)}$.}
    \label{fig:approx_himmelblaur5}
\end{figure}

\paragraph{Himmelblau function.} Here, we consider $\X = \R^2$ and $\Y = \R$, with $A(x)$ being the Himmelblau function. We generate an ensemble with $J=100$ particles according to a Gaussian distribution and pick one arbitrary ensemble member $u^{(i)}$ as a reference point $x$. Figures \ref{fig:approx_himmelblaur1} and \ref{fig:approx_himmelblaur5} show the resulting approximations $A_x$ and $A_x^{(2)}$ in reference to the original function, and to its exact second-order Taylor approximation at $x$. 

\section{Locally weighted ensemble Kalman-type methods}\label{sec:lwEK}
Now that we have set up a solid concept of locally weighted linear approximation, it is very straightforward to define a locally weighted EKI by replacing the empirical covariance in the definition of the linear EK by its weighted variant\new{, see \eqref{eq:locally_weighted_moments} which we provide here again for convenience:\begin{equation*}
\begin{split}
    C_{u,A}^\kappa(x) &= \sum_{j=1}^J \kappa^{(j)}(x) (u^{(j)} - \mu^\kappa(x)) \cdot (A(u^{(j)}) - \mu_A^\kappa(x))^\top\\
    \mu^\kappa(x) &= \sum_{j=1}^J \kappa^{(j)}(x) u^{(j)}\\
    \mu_A^\kappa(x) &= \sum_{j=1}^J \kappa^{(j)}(x) A(u^{(j)})
\end{split}
\end{equation*}
This allows us to define the locally weighted Ensemble Kalman method for inversion:}
\begin{definition} With the definition of locally weighted covariance above, we define the locally weighted EKI (lwEKI) as the system of ODEs given by
      \begin{align*}
 \frac{\d}{\d t} u^{(i)} &= -C_{u,A}^\kappa(u^{(i)}) \cdot \left(A(u^{(i)}) - y\right)
     \end{align*}
 \end{definition}
     This means locally weighted EK just replaces moments like $\mu$, $C_{u,A}$, etc., by their $\kappa$-weighted variants. In a similar way we can modify other algorithms based on the Ensemble Kalman method, such as the EnKF, ALDI, etc.
     
 The hard work in analysing $D_\kappa $ pays off now because the following lemma is trivially true by definition of $D_\kappa $.
\begin{lem}\label{lem:precond_gradflow_locweight}
    The locally weighted EKI can equivalently be defined as   
     \begin{align*}
          \frac{\d}{\d t} u^{(i)} &= -C_{u,u}^\kappa(u^{(i)})  (D_\kappa  A(u^{(i)}))^\top  (A(u^{(i)})-y) \\
          &= - \sum_{j=1}^J \kappa^{(j)}(u^{(i)}) (u^{(j)} - \mu(u^{(i)}))  \left\langle A(u^{(j)}) - \mu_A(u^{(i)}), A(u^{(i)}) - y\right\rangle
     \end{align*}
     \begin{proof}
         \new{Follows by definition and properties of $D_\kappa$.}
     \end{proof}
\end{lem}
Note the similarity of this relation to \eqref{eq:EK_evolution_gradient}. \new{The significance of the result of Lemma \ref{lem:precond_gradflow_locweight} is that this shows that locally weighted EKI is a canonical relaxation of EKI, generalising the property of being a preconditioned gradient flow (which holds for the EKI only in the linear setting) to the nonlinear case, although the flow is following not the true gradient, but its approximation $D_\kappa A$.}

\begin{rem}
    Locally weighted EKI is in general \textit{not} equivalent to the system of ODEs given by
     \begin{align*}
          \frac{\d}{\d t} u^{(i)} &= -C_{u,u}^\kappa(u^{(i)}) D_\kappa  \Phi(u^{(i)}).
     \end{align*}
\end{rem}

It is interesting to see whether these ideas can also be applied to the context of filtering (\cite{reich2015probabilistic,law2015data,chopin2020introduction}). We consider a simplified filtering setting where we focus on the \textit{analysis} step, i.e. the transition from (predicted) prior to (data-updated) posterior. For our purposes, let $\mu_0(\d u)$ be the prior, and $\mu^y(\d u) \propto \exp(-\Phi(u)) \mu_0(\d u)$ be the posterior. A standard homotopy argument (\cite{reich2011dynamical,calvello2022ensemble}) shows that we can view this as a smooth flow of measures via
\[\mu_\tau(\d u) \propto \exp(-\tau\Phi(u)) \mu_0(\d u),\]
with $\tau = 0$ recovering the prior, $\tau = 1$ yielding the posterior $\mu^y$, and $\tau\to \infty$ corresponding to the underregularised EKI inversion setup. Both the Ensemble Square Root filter (EnSRF, \cite{tippett2003ensemble,bergemann2010mollified,bergemann2010localization})
\begin{align*}
    \frac{\d}{\d t} u^{(i)} &= - C_{u,A} \cdot \left( \frac{A(u^{(i)}) + A(m)}{2} - y\right)
\end{align*}
with $m = \frac{1}{J}\sum_{i=1}^J u^{(i)}$; and the stochastic Ensemble Kalman filter (\cite{houtekamer2005ensemble}) 
\begin{align}
    \frac{\d}{\d t} u^{(i)} &= - C_{u,A} \cdot \left(A(u^{(i)})- y + \dot W\right)
\end{align}
can be used in this way to approximate the true posterior $\mu^y$ at time $t=1$, but this will only be valid for linear forward operators $A$. 

By replacing all global moments with their respective locally weighted variants, we obtain an algorithm that can be more robust with respect to nonlinearity of $A$ and multimodality of the measures involved. For example, a locally weighted Ensemble Square Root filter (lwEnSRF) can be defined as
\begin{align}\label{eq:evolution_LWEnSRF}
\frac{\d}{\d t} u^{(i)} &= -C_{u,A}^\kappa(u^{(i)}) \cdot \left(\frac{A(u^{(i)}) + \mu_A^\kappa(u^{(i)})}{2}- y\right)
\end{align}
It should be noted that this evolution needs to be run for a longer time than the non-weighted Ensemble Kalman filter, since the magnitude of the preconditioning factor $C_{u,A}^\kappa$ is systematically smaller than its non-weighted counterpart. This can most easily seen by comparing $C_{u,u}$ and $C_{u,u}^\kappa(x)$ for any point $x$: The locally weighted variant only incorporates physically close particles $u^{(i)}$ in the computation of (locally weighted) covariance, so it will be systematically smaller than the global variance $C_{u,u}$, unless all particles are extremely close together. This means the timescale of the evolution in \eqref{eq:evolution_LWEnSRF} needs to be amplified by a compensating factor.

\subsection{Numerical experiment: locally weighted EKI, 2d}
\begin{figure}
    \centering
    \includegraphics[width=0.5\textwidth]{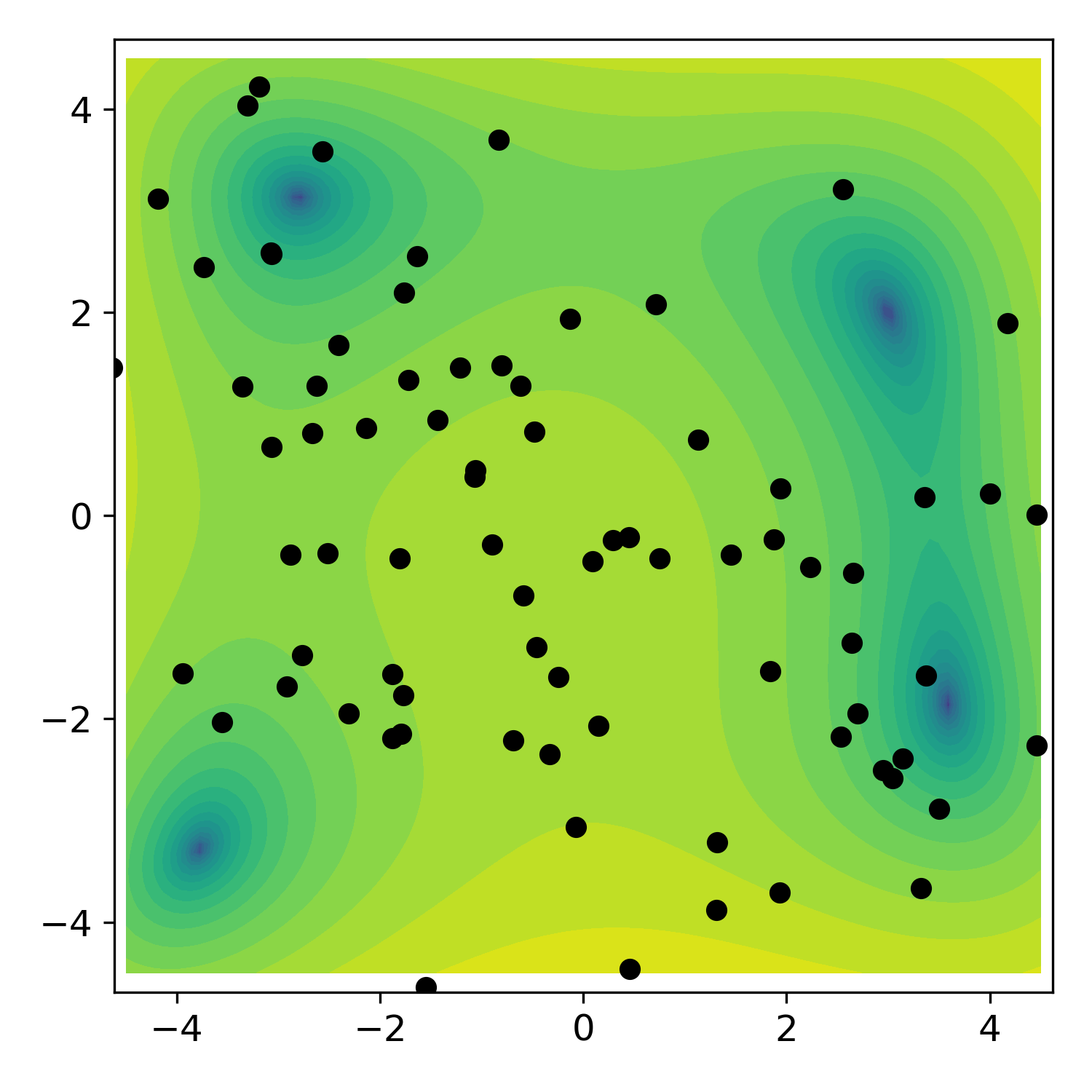}%
    \includegraphics[width=0.5\textwidth]{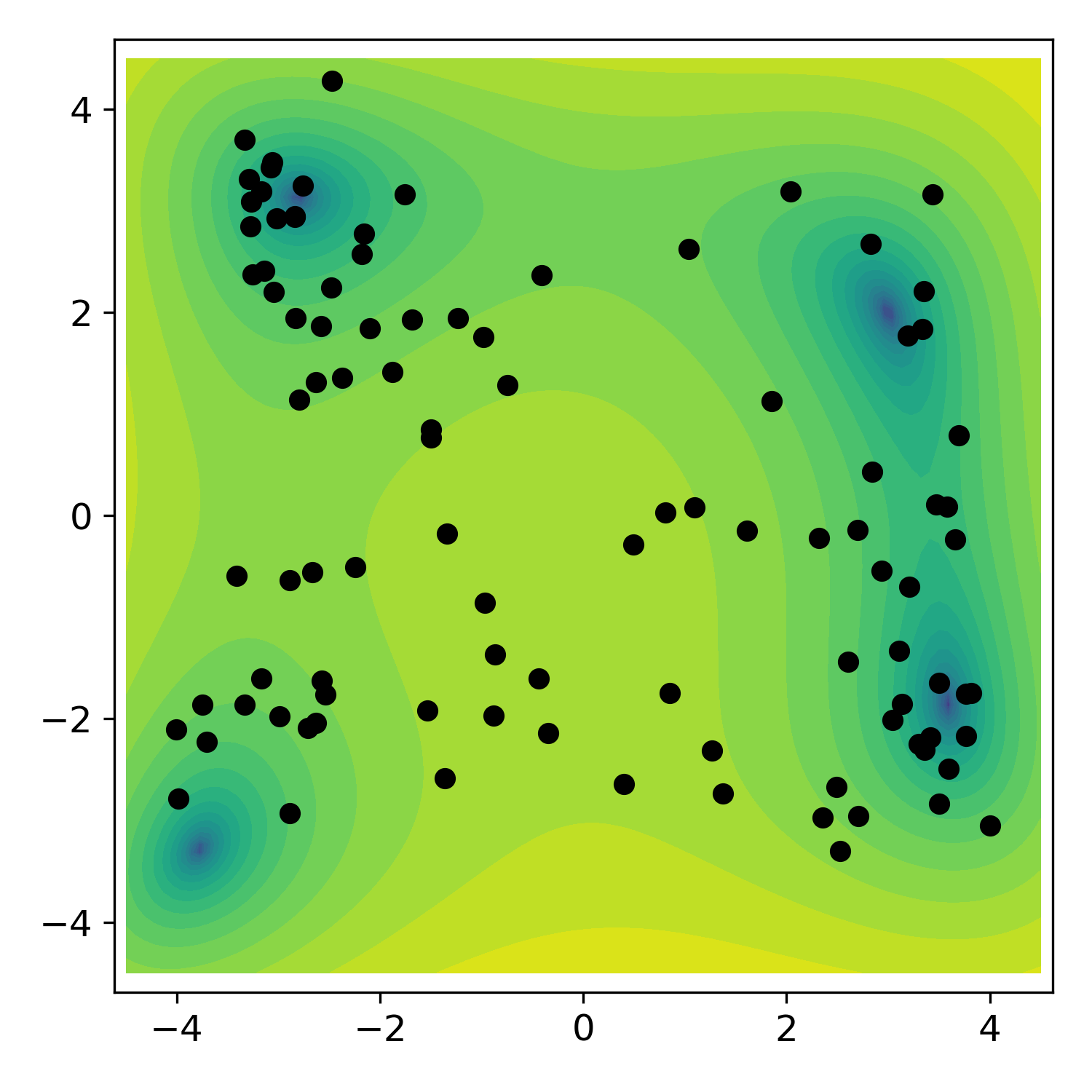}\\
    \includegraphics[width=0.5\textwidth]{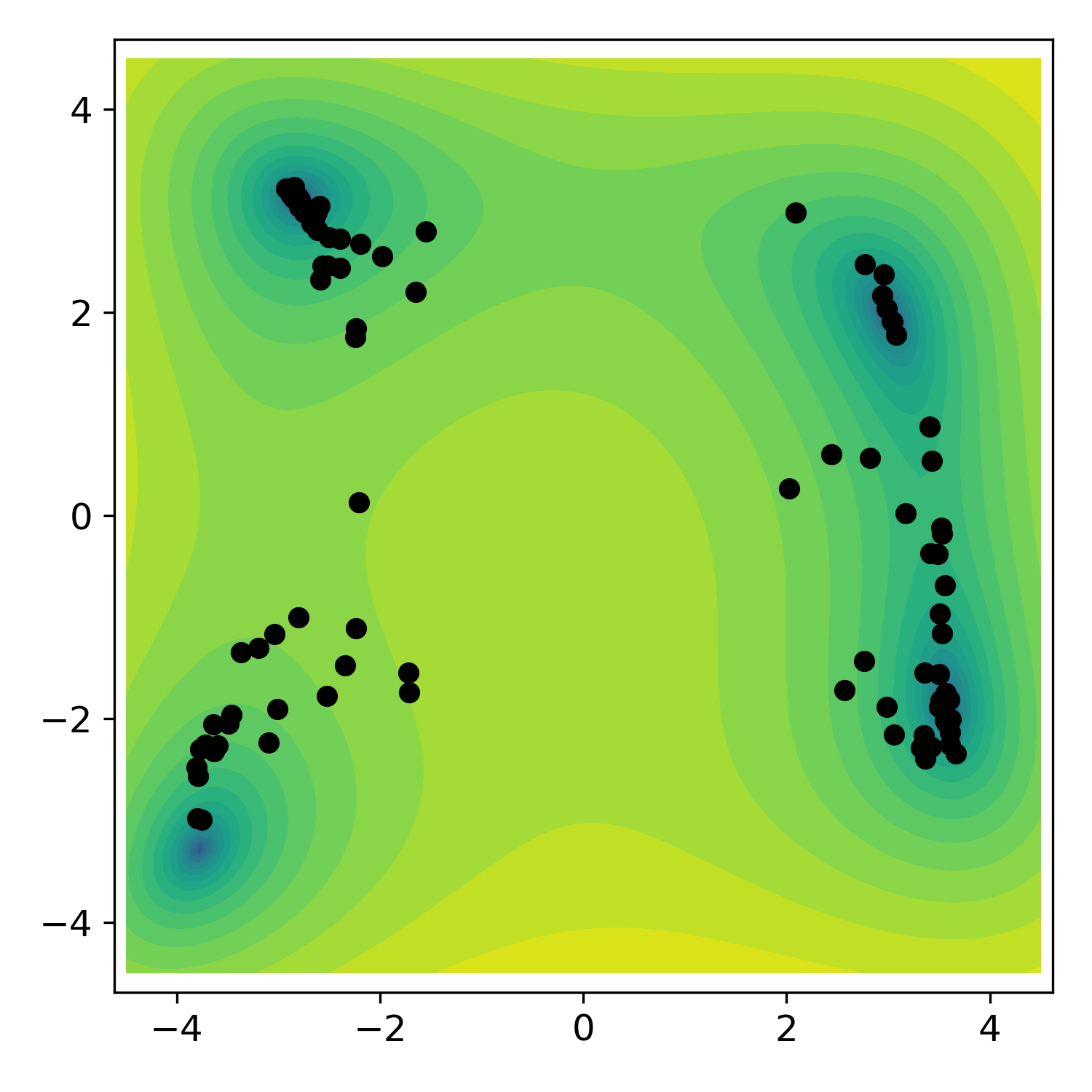}%
    \includegraphics[width=0.5\textwidth]{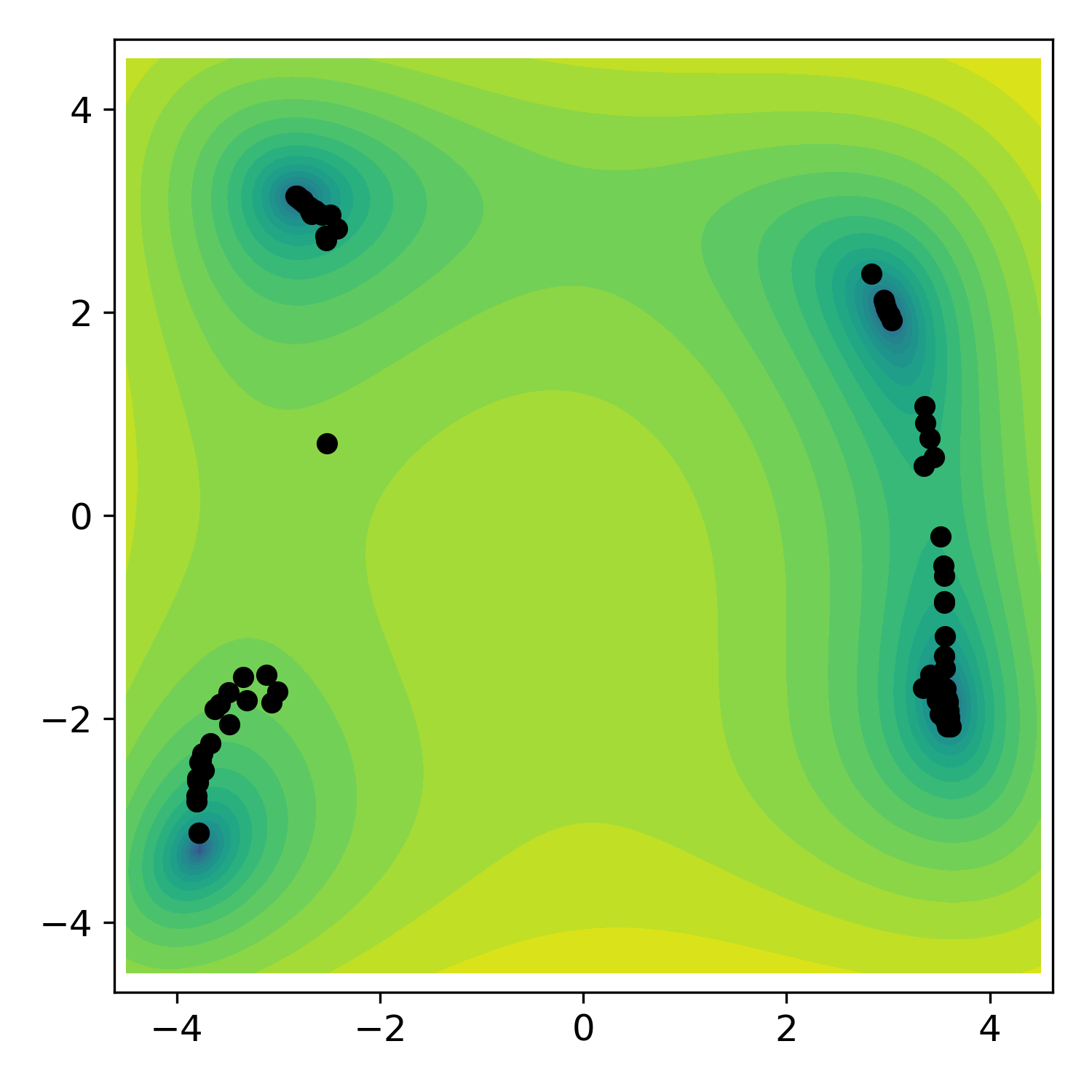}
    \caption{\new{locally weighted EKI for the Himmelblau benchmark problem. Misfit functional $\Phi$ is shown as contour plot, with evolution of ensemble shown for four time steps (in order: top left, top right, bottom left, bottom right). Inversion of $Au = y$ amounts to minimisation of $\Phi$. Ensemble eventually collapses to all four global minima (not shown).}}
    \label{fig:himmelblau_lwEKI}
\end{figure}

We consider a numerical experiment based on the ``Himmelblau'' optimisation benchmark function. Let 
\begin{equation}
    A(x_1,x_2) = (x_1^2 + x_2, x_1 + x_2^2)
\end{equation}
and consider a data/target value $y = (11,7)^\top$. We can define a cost/misfit functional $\Phi(u) = \frac{1}{2}\|y - A(u)\|^2$. Minimisers of this potential are the four solutions to $A(x) = y$ and can be observed as being the minima of the misfit functional $\Phi$ (which happens to be the Himmelblau function), as shown in figure \ref{fig:himmelblau_lwEKI}. We generate an ensemble with $J=100$ particles and perform locally weighted EK with a Gaussian kernel with bandwidth $1$. The results are shown in figure \ref{fig:himmelblau_lwEKI} and it can be observed that the locally weighted EKI successfully finds all four global minima.

\begin{figure}
    \centering
    \includegraphics[width=0.5\textwidth]{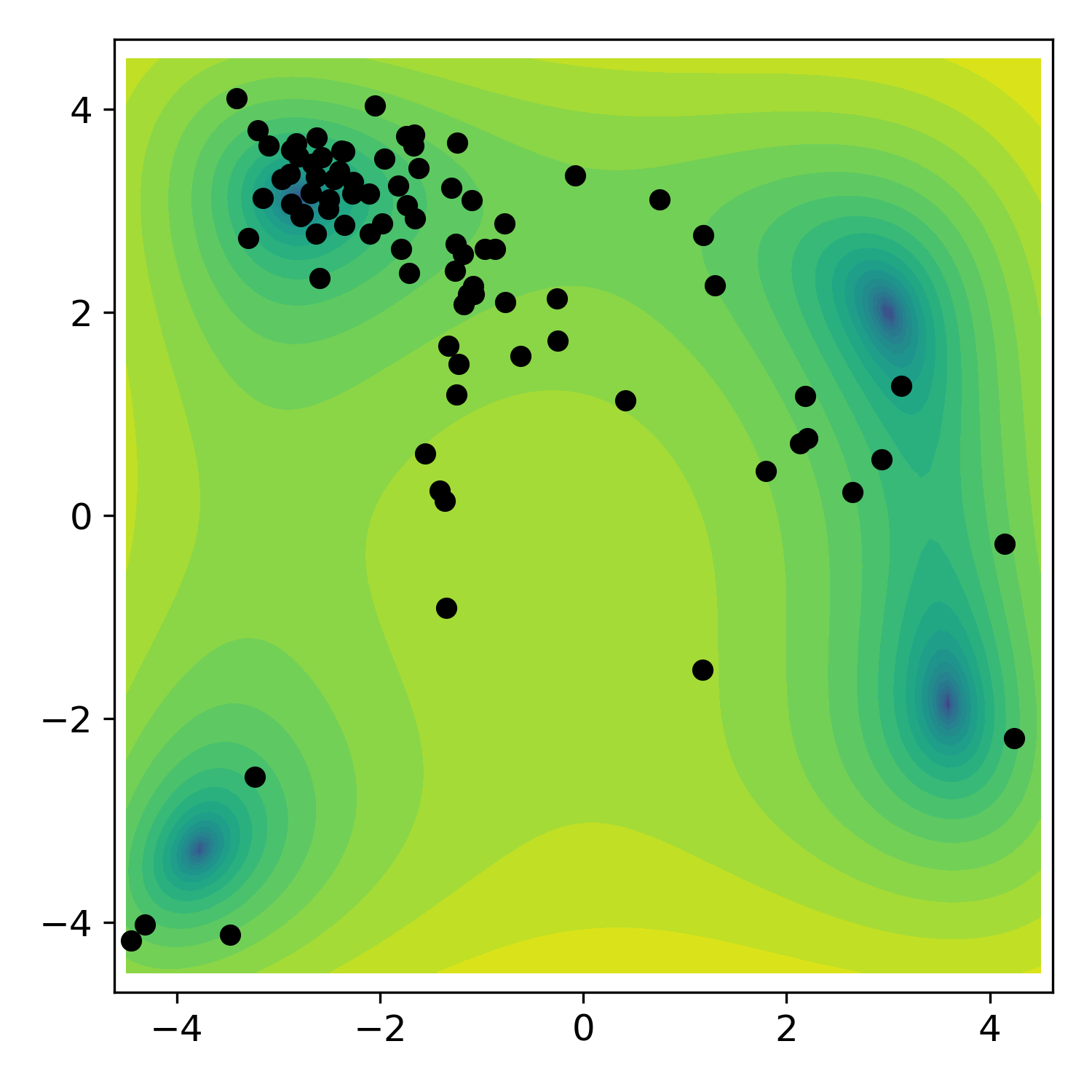}%
    \includegraphics[width=0.5\textwidth]{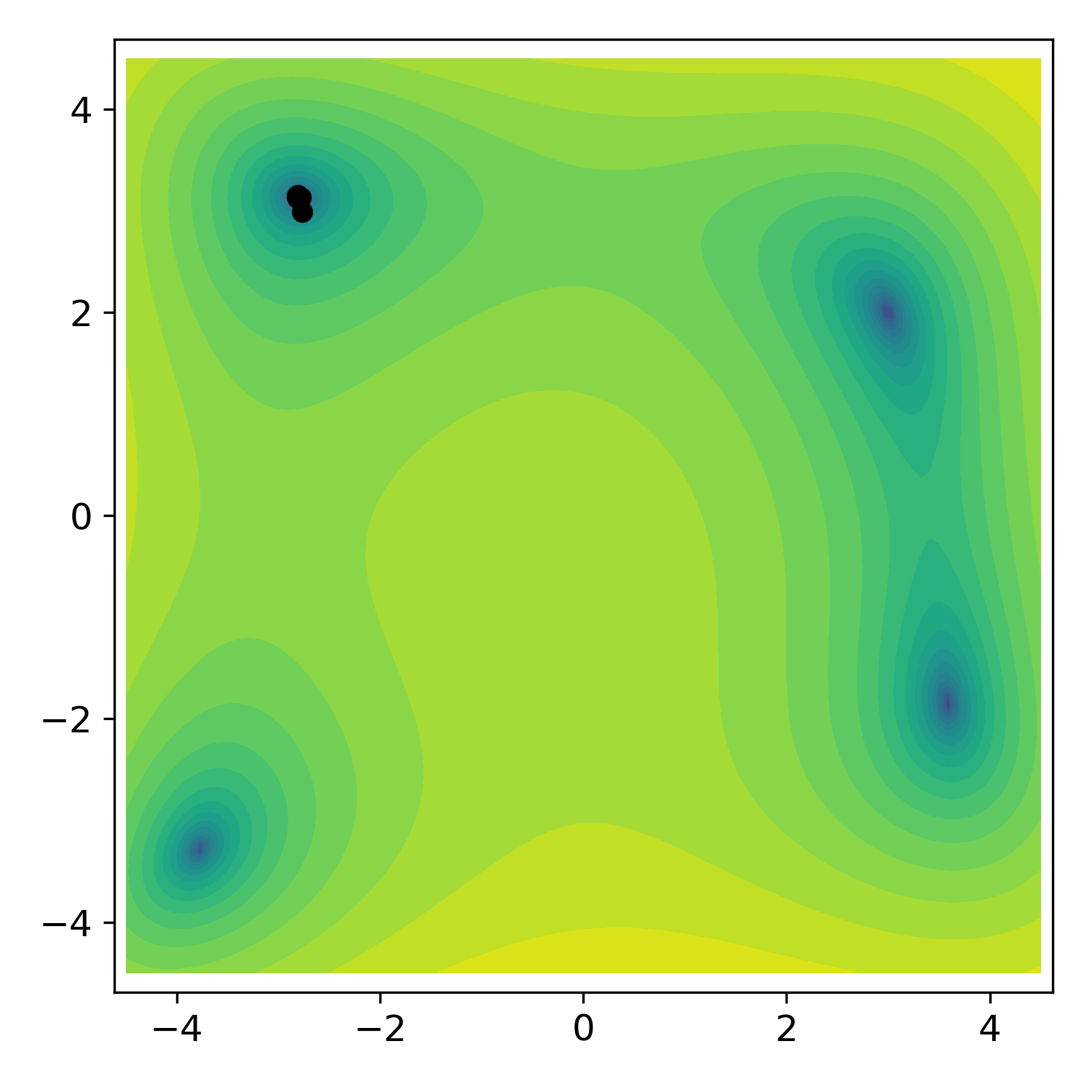}\\
    \caption{\new{non-weighted EKI for the Himmelblau benchmark problem. Shown are two snapshots of (unweighted) EKI. Ensemble eventually collapses to only one of the four global minima.}}
    \label{fig:himmelblau_EKI}
\end{figure}

We compare this with the linear EK in figure \ref{fig:himmelblau_EKI}. It can be observed that the ensemble collapses on one of the four global minima. This should be considered a best-case scenario because it is entirely possible that the linear EKI stops prematurely due to its empirical covariance becoming singular, as discussed above.
\new{This means lwEKI is a promising candidate for optimisation and inversion ontexts, and shows potential to combine the robustness of EKI with a stronger adaptibility to nonlinear and multimodal settings.}

\subsection{Numerical experiment: locally weighted EnSRF, 1d}\label{sec:lwEnSRF}
\begin{figure}
    \centering
    \includegraphics[width=\textwidth]{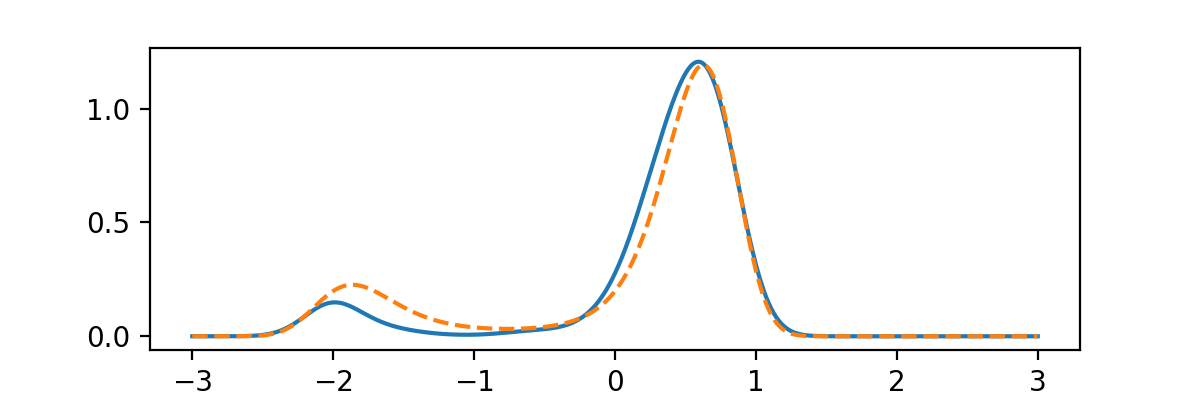}\\
    \includegraphics[width=\textwidth]{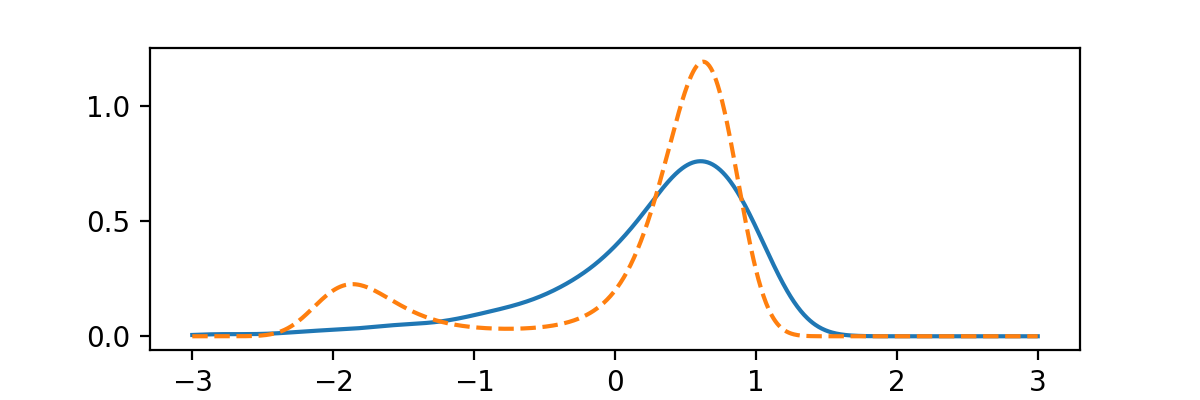}%
    \caption{Side-by-side comparison of (best timestep) posterior approximations of both weighted (top) and unweighted (bottom) EnSRF. True posterior is shown with a dotted line. We \new{draw attention to the fact that the unweighted EnSRF misrepresents both peaks quite seriously.}}
    \label{fig:experiments_sidebyside}
\end{figure}

\begin{figure}
    \centering
    \includegraphics[width=\linewidth]{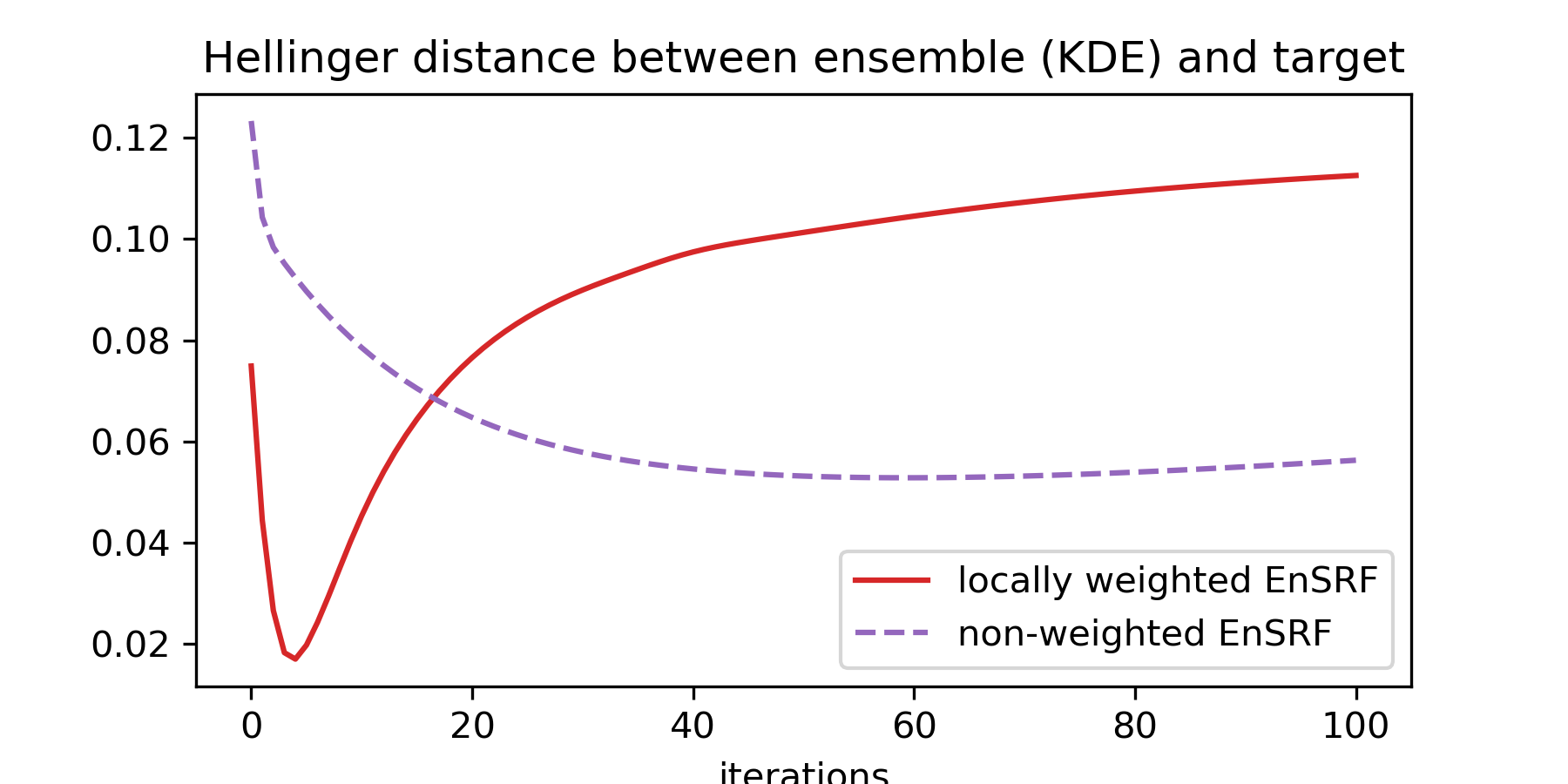}
    \caption{Convergence of Kernel Density Estimator of EnSRF to target distribution, measured in Hellinger distance. We can also see that (no matter whether for locally weighted or non-weighted EnSRF) we need to stop the iteration as soon as suitable distance to target distribution is achieved. }
    \label{fig:Hellinger}
\end{figure}

We consider a simple case similar to the counterexample in \cite{ernst2015analysis}: Let $A(x) = x + 0.75x^2$ and $y = 1$. Even starting with a standard normal Gaussian $N(0,1)$ prior, the posterior will be bimodal. Running the locally weighted EnSRF with an ensemble of 1000 particles and kernel bandwidth $r=0.1$ for $t\in[0,3]$ via an Euler discretisation yields the evolution of measure shown in figure \ref{fig:experiments_sidebyside}, together with the results of an unweighted EnSRF variant. \new{While these results suggest that locally weighted Ensemble Kalman-based methods can be applied for sampling of multimodal distributions, we did experience that some hyperparameter tuning is necessary to make it work, in particular regarding the bandwidth of the kernel used for weighting, and the final time $T$ of the lwENSRF. \neww{Another issue is that the iteration needs to be stopped, to avoid ``overshooting'' of the ensemble, concentrating more than what is required for convergence to the target distribution, see the non-monotonous behaviour in figure \ref{fig:Hellinger}, which shows the Hellinger distance between a kernel density estimation of the ensemble and the target distribution.} In light of this, we recommend the use of the reader's favourite sampling method instead since it is likely to perform better. The benefits of locally weighted ensemble methods seem more apparent in the inversion and optimisation context.}


\subsection{Numerical experiment: locally weighted EKI, 10d}\label{sec:experiment_ring}
For a more moderate-sized example, we consider the symmetrical problem of inferring $u\in\R^{10}$ from noisy norm data $y = \|u\|^2 + \eps$. If we set a centered Gaussian prior $u\sim N(0,\sigma^2 I)$ and $\eps \sim N(0,1)$, then the posterior is given by
\[\mu^y(\d u) \propto \exp\left(-\frac{\|u\|^2}{2\sigma^2} - \frac{(y-\|u\|^2)^2}{2}\right).\]
In order to visually judge the fit of the result of the locally weighted EKI we push the measure forward into data space: If we set $T(u) = \|u\|$, then the distribution of $T(u)|y$ can be seen to be
\[ P(T(u)=dz|y)\propto z^{n-1} e^{-z^2/(2\sigma^2)-(y-z^2)^2/2}. \]

\begin{figure}
    \centering
    \includegraphics[width=0.5\textwidth]{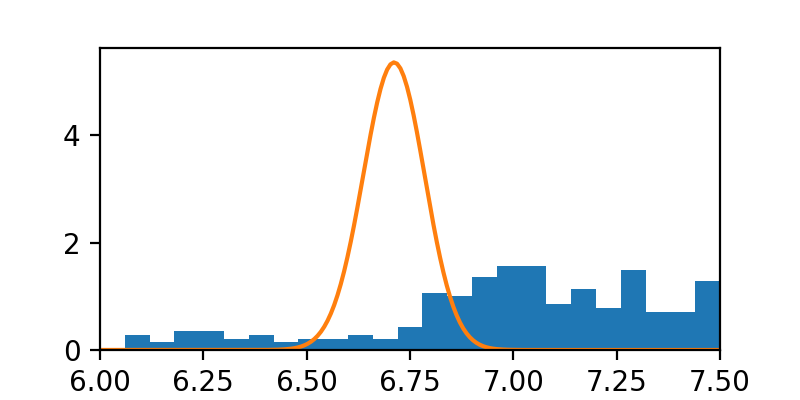}%
    \includegraphics[width=0.5\textwidth]{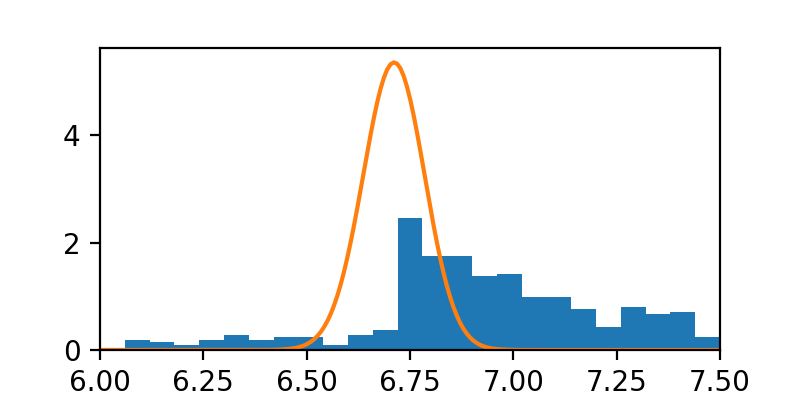}\\
    \includegraphics[width=0.5\textwidth]{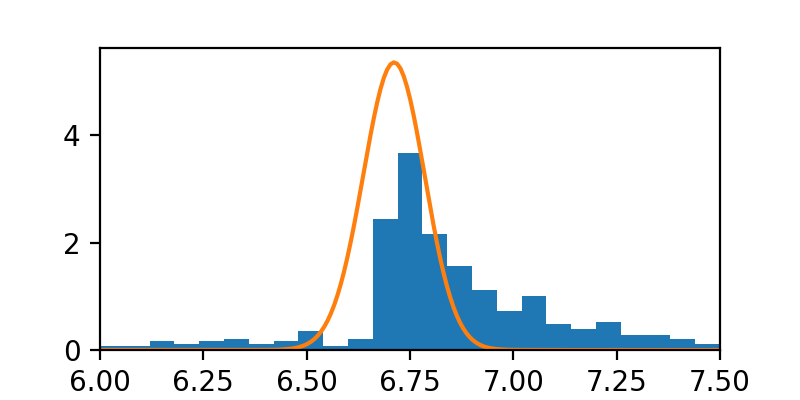}%
    \includegraphics[width=0.5\textwidth]{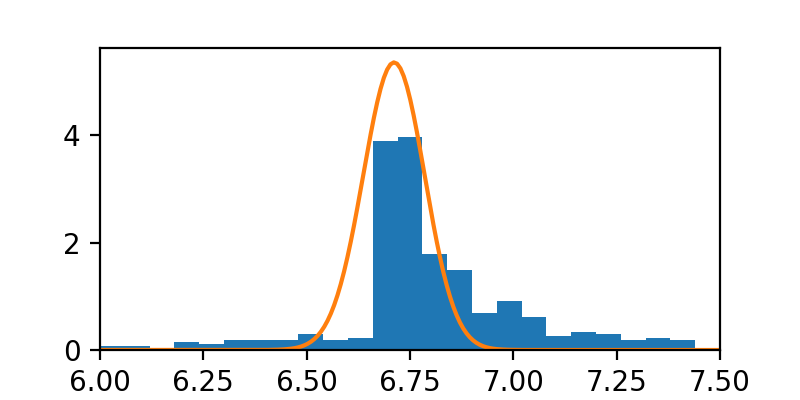}%
    \caption{Ensemble of size $J=500$ for the 10-dimensional problem in section \ref{sec:experiment_ring}, evolving via locally weighted EKI with bandwidth 2. Note that the posterior is not unimodal, but rather concentrated on a thin shell in 10-d, in the vicinity of radius $6.75$ (matching the data $45 = y = \|u\|^2 + \eps$). Figures show histogram Euclidean norm of ensemble members together with theoretical distribution according to posterior measure (orange line). Both particle filter and unweighted Ensemble Kalman completely fail in this example: Prior samples will not hit region of large posterior density at all, which leads to immediate filter degeneracy of the particle filter. The Ensemble Kalman method encounters immediate empirical covariance collapse due to the symmetry of the problem and does not evolve the prior.}
    \label{fig:experiments_ring}
\end{figure}

Note that in this setting both the particle filter and the Ensemble Kalman filter have significant challenges: The effective support of the posterior is a thin shell in 10-dimensional space between radii 6 and 7, and prior samples to lie in this subset is an astronomically rare event, which means the particle filter will immediately lose its complete samples. On the other hand, the forward mapping is such that the empirical covariance $C_{u,A}$ in the initial ensemble will be close to $0$ due to the symmetry of the problem. This means the Ensemble Kalman filter will not move its ensemble at all, failing to recover the posterior. Given that, the locally weighted EKI shows a much better ability to recover an approximation of the posterior, as seen in figure \ref{fig:experiments_ring}, which shows the evolution of a relatively small ensemble of size 500 for a bandwidth of $2$, in the setting where $y = 45$.

\new{
\subsection{Discontinuous optimisation}\label{sec:discontinuous}
Locally weighted EKI can be applied to a particularly challenging setting: Inversion (or optimisation) for discontinuous forward maps. We consider the following inverse problem:
\begin{align*}
    y = \texttt{round}(x)^2 + \varepsilon
\end{align*}
for some Gaussian noise term $\varepsilon$. We set $y = 4$ which means $(-2.5,-1.5)$ and $(1.5,2.5)$ are all minimising the misfit functional $\Phi(x) = \frac{1}{2}\|$. Gradient-based methods do not work here since the misfit is piecewise constant and discontinuous, but the locally weighted EKI is able to represent the solution set, as can be seen in figure \ref{fig:discont}, which shows the result of applying lwEKI with $J=10$ particles to this inverse problem. We can see that the lwEKI evolves the initial ensemble (with no initial particle actually solving the inverse problem already) into the two components $(-2.5,-1.5)$ and $(1.5,2.5)$. This suggests that lwEKI has potential to be applied in optimisation and inversion contexts where the misfit function under consideration is discontinuous or otherwise highly irregular.
\begin{figure}
    \centering
    \includegraphics[width=0.75\linewidth]{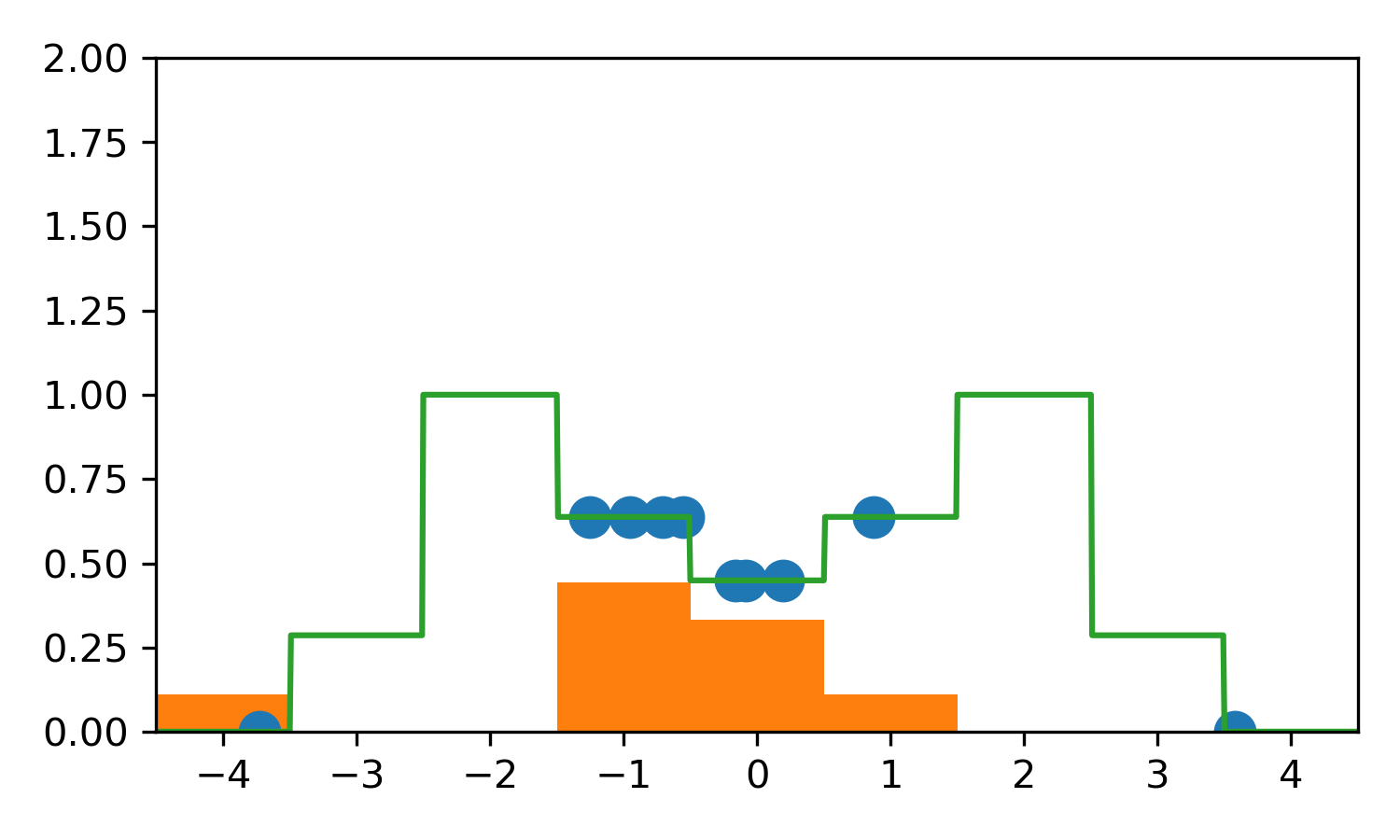}\\
    \includegraphics[width=0.75\linewidth]{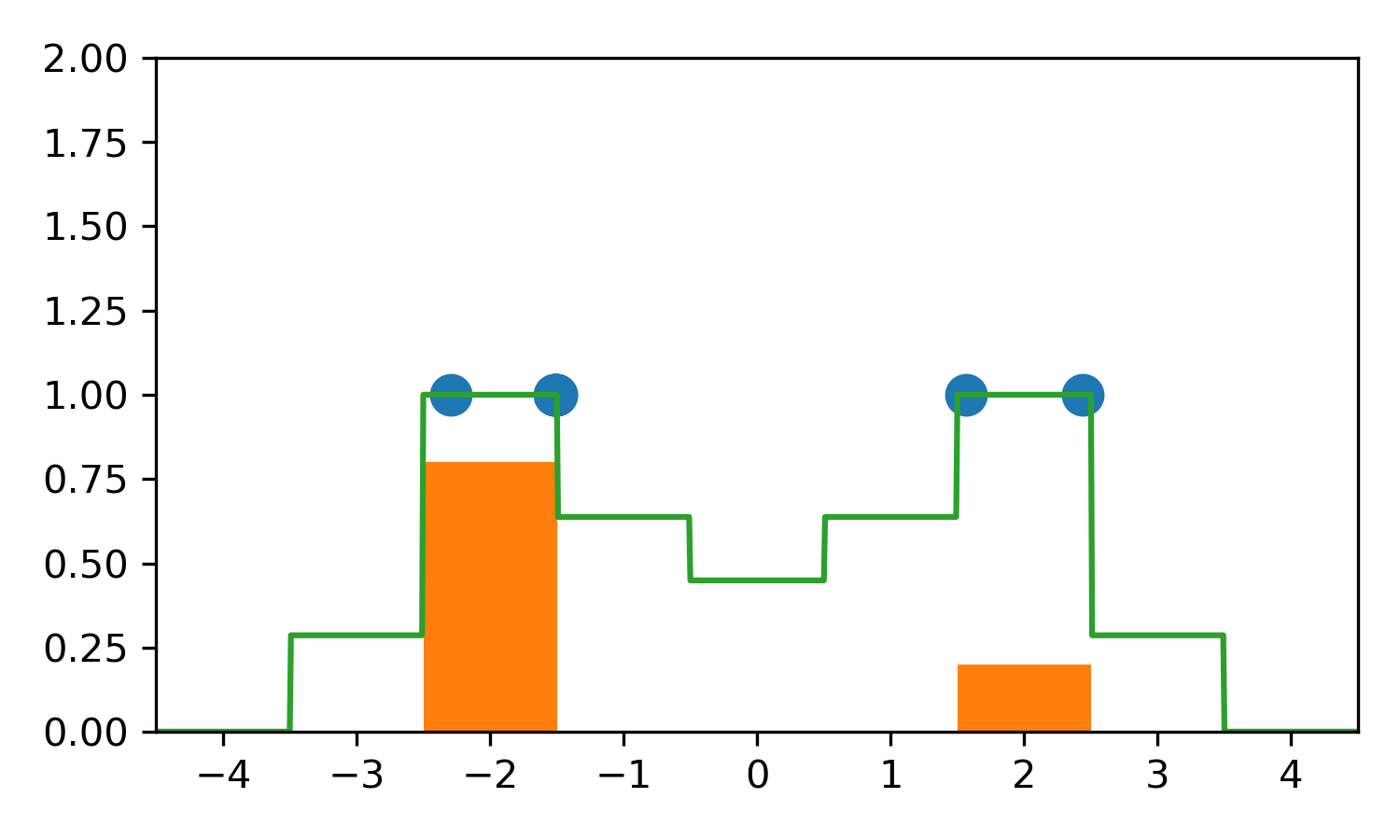}
    \caption{Applying lwEKI to a discontinuous and piecewise constant inverse problem. The green curve shows the likelihood function $\exp(-\Phi)$, and inversion amounts to maximising this function. Top row: Initialisation (ensemble of size 10 is shown both as blue dots, and as a corresponding histogram (orange)). Bottom row: Position of the ensemble after 500 iterations (number of iterations and stepsize were not optimised). Both areas of interest are represented by the ensemble.}
    \label{fig:discont}
\end{figure}}

\new{
\subsection{Computational effort}
We focus on inversion, i.e., the lwEKI algorithm. In analysing computational requirements, there are some issues to note:
\begin{itemize}
    \item The required number of iterations is not clear a priori (similar to other optimisation algorithms) and should be tuned with an automatic stopping criterion.
    \item The dominant computational effort goes into computing the localised covariances. Done naively, this scales like $\mathcal O(J^2)$ (since every particle's covariance is computed by computing a weighted average over all other particles). This can be mitigated by computation tricks (see a discussion in the conclusion), but remains the crux for this method. 
\end{itemize}
}

\subsection{Conclusion, \new{computational} challenges and outlook}
The elephant in the room regarding the locally weighted EK is the larger computational overhead (as compared to the linear EK): Since moments have to be generated for each particle, rather than just for the whole ensemble, this creates considerable computational cost for large ensemble sizes. There are several ways of mitigating this issue: For a sufficiently quickly decaying kernel, the locally weighted covariances have a very low-rank structure, which can possibly be used for computational speed-up. Alternatively, we can further modify the locally weighted EKI using ideas of the ``cluster-based'' variant of polarCBO in \cite{bungert2022polarized}: By fixing a number of clusters $M$, and relating all particles to clusters instead of all particles with each other, we only need $\mathcal O(M\cdot J)$ comparisons, rather than $\mathcal O(J^2)$. This idea was already proposed by \cite{wagner2022ensemble} in the context of rare event estimation. Thirdly, we can use the kernel to build a graph structure on the ensemble, where there is an edge between particles $u^{(i)},u^{(j)}$ if the weight $k(u^{(i)},u^{(j)})$ is positive (or larger than some threshold). Then graph-based message-passing algorithms can be utilised to improve efficiency of computation. \new{A third option to reduce computational complexity of the covariance computation might be ideas similar to the extrapolation idea in \cite{schillings2023ensemble}: Instead of estimating the covariance locally everywhere, we estimate it at only selected particles, and extrapolate this information (using a Taylor approximation) to particles in the vicinity.}

But not only large value of ensemble size $J$ are problematic: If the dimensionality of parameter space is very high, radial distance-based local weighting will have performance issues due to the curse of dimensionality: Since relative volumes of balls of fixed radius have asymptotically vanishing values, it becomes harder and harder to find ensemble members which are not getting weighted down to $0$ by the kernel function. This can be remedied by, e.g., utilising lower-dimensional structure or replacing distance-based local weights by something more robust.

Subject of future work is a more realistic feasibility test of the locally weighted Ensemble Kalman method in a realistic filter setting. It would be interesting to see whether the locally weighted EnKF is able to combine the speed of the EnKF with the flexibility of a particle filter in a nonlinear filtering setup. Open questions pertain the stopping time for algorithms like the lwEnSRF (since $T=1$ is not tied to the posterior measure anymore, for reasons explained above), and more generally bias reduction methods. Possible extensions include the use of locally weighted Ensemble Kalman methods as proposal measures to be corrected with a Metropolis acceptance step, or by reweighting using the correct likelihood.

Another interesting question is whether the kernel bandwidth should be kept fixed throughout, adapted, or even locally adapted (for example based on the locally weighted covariance, closing this loop, too).
We refer further elaboration of these considerations to future research.

\section*{Acknowledgments}
Some ideas in this manuscript were inspired by conversations with Matei Hanu, Claudia Schillings, Björn Sprungk, Chris Stevens, and Simon Weißmann. The author acknowledges support from MATH+ project EF1-19: Machine Learning Enhanced Filtering Methods for Inverse Problems, funded by the Deutsche Forschungsgemeinschaft (DFG, German Research Foundation) under Germany’s Excellence Strategy – The Berlin Mathematics Research Center MATH+ (EXC-2046/1, project ID: 390685689).

\section*{Statements and Declarations}
Competing Interests: Nothing to declare.
\printbibliography

\end{document}